\documentclass{elsart}

\usepackage{amssymb}
\newtheorem{lemma}{Lemma}[section]   
\newtheorem{theorem}{Theorem}[section]    \newtheorem{proposition}{Proposition}[section]
\newtheorem{definition}{Definition}[section]     \newtheorem{remark}{Remark}[section]
\def\1{\l\!\l}

\begin{document}

\begin{frontmatter}



\title{A Mixed Generalized Multifractal Formalism For Vector Valued Measures}
\author{Anouar Ben Mabrouk}
\address{Computational Mathematics Laboratory, Department of Mathematics\\ Faculty of Sciences, 5019 Monastir, Tunisia.}
\ead{anouar.benmabrouk@issatso.rnu.tn}
\begin{abstract}
We introduce a mixed generalized multifractal formalism which extends the mixed multifractal formalism introduced by L. Olsen based on generalizations of the Hausdorff and packing measures. The validity of such a formalism is proved in some special cases.
\end{abstract}
\end{frontmatter}
\section{Introduction and main results}

$Dom(B)=\{-\nabla\,B_\mu(q);\,\nabla\,B_\mu\,\exists\,\}$ and $f_\mu$

\section{Hausdorff and packing measures and dimensions}
Given a subset $E\subseteq\mathbb{R}$, and $\epsilon>0$, we call an $\epsilon$-covering of $E$, any countable set $(U_i)_i$ of non-empty subsets $U_i\subseteq\mathbb{R}$ satisfying
\begin{equation}\label{epsiloncovering}
E\subseteq\displaystyle\bigcup_iU_i\quad\hbox{and}\quad|U_i|=diam(U_i)\leq\epsilon,
\end{equation}
where for any subset $U\subseteq\mathbb{R}$, $|U|=diam(U)$ is the diameter defined by $|U|=diam(U)=\displaystyle\sup_{x,y\in\,U}|x-y|$. Remark here that for $\epsilon_1<\epsilon_2$, any $\epsilon_2$-covering of $E$ is obviously an $\epsilon_2$-covering of $E$. This implies that the quantity
$$
\mathcal{H}^s_\epsilon(E)=\inf\{\displaystyle\sum_i|U_i|^s\,\,;\,\,(U_i)\,\,\hbox{satisfying}\,\,\ref{epsiloncovering}\}
$$
is a non increasing function in $\epsilon$. It's limit
$$
\mathcal{H}^s(E)=\lim_{\epsilon\downarrow0}\mathcal{H}^s_\epsilon(E)
$$
defines the so-called $s$-dimensional Hausdorff measure of $E$. It holds that for any set $E\subseteq\mathbb{R}$ there exists a critical value $s_E$ in the sense that
$$
\mathcal{H}^s(E)=0,\,\,\forall\,s<s_E\quad\hbox{and}\quad\mathcal{H}^s(E)=+\infty,\,\,\forall\,s>s_E,
$$
or otherwise,
$$
s_E=\displaystyle\sup\{s>0\,;\,\,\mathcal{H}^s(E)=0\}=\displaystyle\inf\{s>0\,;\,\,\mathcal{H}^s(E)=+\infty\}.
$$
Such a value is called the Hausdorff dimension of the set $E$ and is usually denoted by $dim_H\,E$ or simply $dim\,E$. When $U_i=B(x_i,r_i)$ is a ball centered at $x_i\in\,E$ and with diameter $r_i<\epsilon$, the covering $(B(x_i,r_i))_i$ is called an $\epsilon$-centered covering of $E$. However, surprisingly, the quantity $\mathcal{H}^s$ restricted only on centered coverings does not define a measure. To obtain a good measure with centered coverings on should do more. Denote
$$
\overline{\mathcal{C}}^s_\epsilon(E)=\inf\{\displaystyle\sum_i|2r_i|^s\,\,;\,\,(B(x_i,r_i))_i\,\,\hbox{an}\,\epsilon-\hbox{centered covering of}\,E\}
$$
and similarly as above,
$$
\overline{\mathcal{C}}^s(E)=\lim_{\epsilon\downarrow0}\overline{\mathcal{C}}^s_\epsilon(E).
$$
As stated previously, this is not a good measure. Indeed, .............

So, to obtain a good candidate, we set for $E\subseteq\mathbb{R}$,
$$
\mathcal{C}^s(E)=\displaystyle\sup_{F\subseteq\,E}\overline{\mathcal{C}}^s(F).
$$
It is called the centered Hausdorff $s$-dimensional measure of $E$. But, although a fascinating relation to the Hausdorff measure exists. It holds that
\begin{equation}\label{RecCentreEtHaudorffMeas}
2^{-s}\mathcal{C}^s(E)\leq\mathcal{H}^s(E)\leq\mathcal{C}^s(E);\,\,\forall\,E\subseteq\mathbb{R}^d.
\end{equation}
Indeed, let $F\subseteq\,E$ be subsets of $\mathbb{R}^d$. It follows from the definition of $\mathcal{H}^s$ and $\overline{\mathcal{C}}^s$ that  $\mathcal{H}^s(F)\leq\overline{\mathcal{C}}^s(F)$. Next, from the fact that $\mathcal{H}^s$ is an outer metric measure on $\mathbb{R}^d$, and the definition of $\mathcal{C}^s$, il results that $\mathcal{H}^s(E)\leq\mathcal{C}^s(E)$. Next, let $\{U_j\}_j$ be an $\epsilon$-covering of $F$ and $r_j=diam(U_j)$. For each $i$ fixed, consider a point $x_i\in\,U_i\cap\,F$. This results in a centered $\epsilon$-covering $\{B(x_i,r_i)\}_i$ of $F$. Consequently,
$$
\overline{\mathcal{C}}_\epsilon^s(F)\leq\displaystyle\sum_i(2r_i)^s=2^s\displaystyle\sum_i(diam(U_i))^s.
$$
Hence,
$$
\overline{\mathcal{C}}_\epsilon^s(F)\leq2^s\mathcal{H}_\epsilon^s(F).
$$
Next, as $\epsilon\downarrow0$, we obtain
$$
\overline{\mathcal{C}}^s(F)\leq2^s\mathcal{H}^s(F),\,\,\forall\,F\subseteq\,E.
$$
which guaranties that
$$
\mathcal{C}^s(E)\leq2^s\mathcal{H}^s(E).
$$
It holds that these measures give rise to some critical values in the sense that, for any set $E\subseteq\mathbb{R}$ there exists a critical value $h_E$ and $c_E$ for which
$$
\mathcal{H}^s(E)=0,\,\,\forall\,s<h_E\quad\hbox{and}\quad\mathcal{H}^s(E)=+\infty,\,\,\forall\,s>h_E
$$
and similarly
$$
\mathcal{C}^s(E)=0,\,\,\forall\,s<c_E\quad\hbox{and}\quad\mathcal{C}^s(E)=+\infty,\,\,\forall\,s>c_E.
$$
But using equation \ref{RecCentreEtHaudorffMeas} above, it proved that $h_E=c_E$ and otherwise,
$$
h_E=\displaystyle\sup\{s>0\,;\,\,\mathcal{H}^s(E)=0\}=\displaystyle\inf\{s>0\,;\,\,\mathcal{H}^s(E)=+\infty\}.
$$
Such a value is called the Hausdorff dimension of the set $E$ and is usually denoted by $dim_H\,E$ or simply $dim\,E$.

Similarly, we call a centered $\epsilon$-packing of $E\subseteq\mathbb{R}^d$, any countable set $(B(x_i,r_i))_i$ of disjoint balls centered at points $x_i\in\,E$ and with diameters $r_i<\epsilon$. The packing measure and dimension are defined as
follows.
$$
\displaystyle\overline\mathcal{P}^s(E)=\lim_{\varepsilon\downarrow0}\Bigl(\sup\{\displaystyle\sum_i(2r_i)^s\,;\,\,(B(x_i,r_i))_i\,\epsilon-\hbox{packing of}\,E\}\Bigr),
$$
$$
\displaystyle{\cal P}^s(E)=\inf\{\displaystyle\sum_i\overline\mathcal{P}^s(E_i)\,\,;\,E\subseteq\cup_iE_i\}.
$$
It holds as for the Hausdorff measure that there exists critical values $\Delta_E$ and $p_E$ satisfying respectively
$$
{\overline\mathcal{P}}^s(E)=+\infty\,\hbox{for}\,s<\Delta(E)\quad\hbox{and}\quad{\overline\mathcal{P}}^s(E)=0\,\hbox{for}\,\alpha>\Delta(E)
$$
and respectively
$$
\mathcal{P}^s(E)=\infty\,\hbox{for}\,s<p_E\quad\hbox{and}\quad\mathcal{P}^s(E)=0\,\hbox{for}\,s>p_E.
$$
The critical value $\Delta(E)$ is called the logarithmic index of $E$ and $p_E$ is called the packing dimension of $E$ denote by $Dim_P(E)$ or simply $Dim(E)$. These quantities may be shown as
$$
\Delta(E)=\sup\{s\,;\,\,\overline{\mathcal{P}}^s(E)=0\}=\inf\{s\,;\,\,\overline{\mathcal{P}}^s(E)=+\infty\}.
$$
and respectively
$$
Dim(E)=\sup\{s\,;\,\,\mathcal{P}^s(E)=0\}=\inf\{s\,;\,\,\mathcal{P}^s(E)=+\infty\}.
$$
Usually, we have the inequality
$$
dim(E)\leq\,Dim(E)\leq\Delta(E),\,\,\forall\,E\subseteq\mathbb{R}^d.
$$
\begin{definition}\label{TaylorDefFractalSets}
A set $E\subseteq\mathbb{R}^d$ is said to be fractal in the sense of Taylor iff $dim(E)=Dim(E)$.
\end{definition}
\section{Multifractal generalizations of Hausdorff and packing measures}\label{ghpmd}
Let $\mu$ be a Borel probability measure on $\mathbb{R}^d$, and a nonempty set $E\subseteq\mathbb{R}^d$ and $\epsilon>0$. Let also $q,\,t$ be real numbers. We will recall hereafter the steps leading to the multifractal generalizations of the Hausdorff and packing measures due to L. olsen in \cite{olsen1}. Denote
$$
{\overline\mathcal{H}}_{\mu,\epsilon}^{q,t}(E)=\displaystyle\inf\{\,\sum_i(\mu(B(x_i,r_i)))^q(2r_i)^t\,\},
$$
where the inf is taken over the set of all centered $\epsilon$-coverings of $E$, and for the empty set, ${\overline\mathcal{H}}_{\mu,\epsilon}^{q,t}(\emptyset)=0$. As for the preceding cases of Hausdorff and packing measures, it consists of a non increasing quantity as a function of $\varepsilon$. We then consider its limit
$$
{\overline\mathcal{H}}_\mu^{q,t}(E)=\displaystyle\lim_{\epsilon\downarrow0}{\overline\mathcal{H}}_{\mu,\epsilon}^{q,t}(E)
=\sup_{\delta>0}{\overline\mathcal{H}}_{\mu,\epsilon}^{q,t}(E)
$$
and finally, the multifractal generalization of the $s$-dimensional Huasdorrf measure
$$
\mathcal{H}_\mu^{q,t}(E)=\displaystyle\sup_{F\subseteq\,E}{\overline\mathcal{H}}_\mu^{q,t}(F).
$$
Similarly, we define the multifractal generalization of the packing measure as follows.
$$
{\overline\mathcal{P}}_{\mu,\epsilon}^{q,t}(E)=\displaystyle\sup\{\,\sum_i(\mu(B(x_i,r_i)))^q(2r_i)^t\,\}
$$
where the sup is taken over the set of all centered $\epsilon$-packings of $E$. For the empty set, we set as usual ${\overline\mathcal{P}}_{\mu,\epsilon}^{q,t}(\emptyset)=0$. Next,
$$
{\overline\mathcal{P}}_\mu^{q,t}(E)=\displaystyle\lim_{\epsilon\downarrow0}{\overline\mathcal{P}}_{\mu,\epsilon}^{q,t}(E)
=\inf_{\delta>0}{\overline\mathcal{P}}_{\mu,\epsilon}^{q,t}(E)
$$
and finally,
$$
\mathcal{P}_\mu^{q,t}(E)=\displaystyle\inf_{E\subseteq\,\cup_iE_i}\sum_i{\overline\mathcal{P}}_\mu^{q,t}(E_i).
$$
In \cite{olsen1}, it has been proved that the measures $\mathcal{H}_{\mu}^{q,t}$, $\mathcal{P}_\mu^{q,t}$ and the pre-measure ${\overline\mathcal{P}}_\mu^{q,t}$ assign in a usual way a dimension to every set $E\subseteq\mathbb{R}^d$ as resumed in the following proposition.
\begin{proposition}\cite{olsen1}
Given a subset $E\subseteq\mathbb{R}^d$,
\begin{enumerate}
\item There exists a unique number $dim_\mu^q(E)\in[-\infty,+\infty]$ such that
$$
\mathcal{H}_\mu^{q,t}(E)=\left\{\matrix{+\infty&\hbox{for}&t<dim_\mu^q(E)\hfill\cr\medskip
0&\hbox{si}&t>dim_\mu^q(E)\hfill}\right.
$$
\item There exists a unique number $Dim_\mu^q(E)\in[-\infty,+\infty]$ such that
$$
\mathcal{P}_\mu^{q,t}(E)=\left\{\matrix{+\infty&\hbox{for}&t<Dim_\mu^q(E)\hfill\cr\medskip
0&\hbox{for}&t>Dim_\mu^q(E)\hfill}\right.
$$
\item There exists a unique number $\Delta_\mu^q(E)\in[-\infty,+\infty]$ such that
$$
{\overline\mathcal{P}}_\mu^{q,t}(E)=\left\{\matrix{+\infty&\hbox{for}&t<\Delta_\mu^q(E)\hfill\cr\medskip
0&\hbox{for}&t>\Delta_\mu^q(E)\hfill}\right.
$$
\end{enumerate}
\end{proposition}
The quantities $dim_\mu^q(E)$, $Dim_\mu^q(E)$ and $\Delta_\mu^q(E)$ defines the so-called multifractal generalizations of the Hausdorff dimension, the packing dimension and the logarithmic index of the set $E$. More precisely, one has
$$
dim_\mu^0(E)=dim(E),\quad\,Dim_\mu^0(E)=Dim(E)\quad\hbox{and}\quad\Delta_\mu^0(E)=\Delta(E).
$$
The characteristics of these functions have been studied completely by L. Olsen. He proved among author results that $dim_\mu^q$ and $Dim_\mu^q$ are monotones and $\sigma$-stables. Furthermore, if $E=support(\mu)$ is the support of the measure $\mu$, one obtains
\begin{description}
\item[a.] The functions $q\longmapsto\,Dim_{\mu}^q(E)$ and $q\longmapsto\Delta_{\mu}^q(E)$ are convex non increasing.
\item[b.] $q\longmapsto\,dim_{\mu}^q(E)$ is non increasing.
\item[c.]
\begin{description}
  \item[i.] For $q<1$; $0\leq\,dim_\mu^q(E)\leq\,Dim_\mu^q(E)\leq\Delta_\mu^q(E)$.
  \item[ii.] $dim_\mu^1(E)\leq\,Dim_\mu^1(E)\leq\Delta_\mu^1(E)=0$.
  \item[iii.] For $q>1$; $dim_\mu^q(E)\leq\,Dim_\mu^q(E)\leq\Delta_\mu^q(E)\leq0$.
\end{description}
\end{description}
\section{Mixed multifractal generalizations of Hausdorff and packing measures and dimensions}
The purpose of this section is to present our ideas. As it is noticed from the literature on multifractal analysis of measures, this latter always considered a single measure and studies it scaling behavior as well as the multifractal formalism associated. Recently, many works have been focused on the study of simultaneous behaviors of finitely many measures. In \cite{olsen3}, a mixed multifractal analysis is developed dealing with a generalization of R\'enyi dimensions for finitely many self similar measures. This was one of the motivations leading to our present paper. Secondly, we intend to combine the generalized Hausdorff and packing measures and dimensions recalled in section \ref{ghpmd} with Olsen's results in \cite{olsen3} to define and develop a more general multifractal analysis for finitely many measures by studying their simultaneous regularity, spectrum and to define a mixed multifractal formalism which may describe better the geometry of the singularities's sets of these measures especially simultaneous singularities.\\
Let $\mu=(\mu_1,\,\mu_2,\dots,\mu_k)$ a vector valued measure composed of probability measures on $\mathbb{R}^d$. We aim to study the simultaneous scaling behavior of $\mu$ which we denote
$$
\displaystyle\lim_{r\downarrow0}\frac{\log\mu(B(x,r))}{\log\,r}\equiv\bigl(\displaystyle\lim_{r\downarrow0}\frac{\log\mu_1(B(x,r))}{\log\,r},
\dots,\,\displaystyle\lim_{r\downarrow0}\frac{\log\mu_k(B(x,r))}{\log\,r}\bigr).
$$
In this paper, we apply the techniques of L. Olsen especially in \cite{olsen1} and \cite{olsen3} with the necessary modifications to give a detailed study of computing general mixed multifractal dimensions of simultaneously many finite number of measures and try to project our results for the case of a single measure to show the generecity of our's. Let a $E\subseteq\mathbb{R}^d$ be a nonempty set and $\epsilon>0$. Let also $q=(q_1,q_2,\dots,q_k)\in\mathbb{R}^k$ and $t\in\mathbb{R}$. The mixed generalized multifractal Hausdorff measure is defined as follows. Denote
$$
\mu(B(x,r))\equiv\bigl(\mu_1(B(x,r)),\dots,\,\mu_k(B(x,r))\bigr)
$$
and the product
$$
(\mu(B(x,r)))^q\equiv(\mu_1(B(x,r)))^{q_1}\dots(\mu_k(B(x,r)))^{q_k}.
$$
Denote next,
$$
{\overline\mathcal{H}}_{\mu,\epsilon}^{q,t}(E)=\displaystyle\inf\{\,\sum_i(\mu(B(x_i,r_i)))^q(2r_i)^t\,\},
$$
where the inf is taken over the set of all centered $\epsilon$-coverings of $E$, and for the empty set, ${\overline\mathcal{H}}_{\mu,\epsilon}^{q,t}(\emptyset)=0$. As for the single case, of Hausdorff measure, it consists of a non increasing function of the variable $\varepsilon$. So that, its limit as $\epsilon\downarrow0$ exists. Let
$$
{\overline\mathcal{H}}_\mu^{q,t}(E)=\displaystyle\lim_{\epsilon\downarrow0}{\overline\mathcal{H}}_{\mu,\epsilon}^{q,t}(E)
=\sup_{\delta>0}{\overline\mathcal{H}}_{\mu,\epsilon}^{q,t}(E).
$$
Let finally
$$
\mathcal{H}_\mu^{q,t}(E)=\displaystyle\sup_{F\subseteq\,E}{\overline\mathcal{H}}_\mu^{q,t}(F).
$$
\begin{lemma}\label{anouarlemme1}
$\mathcal{H}_\mu^{q,t}$ is an outer metric measure on $\mathbb{R}^d$.
\end{lemma}
\ {\bf Proof.} We will prove firstly that $\mathcal{H}_\mu^{q,t}$ is an outer measure. This means that
\begin{description}
\item [i.] $\mathcal{H}_\mu^{q,t}(\emptyset)=0$.
\item[ii.] $\mathcal{H}_\mu^{q,t}$ is monotone, i.e. $\mathcal{H}_\mu^{q,t}(E)\leq\mathcal{H}_\mu^{q,t}(F)$, whenever $E\subseteq\,F\subseteq\mathbb{R}^d$.
\item[iii.] $\mathcal{H}_\mu^{q,t}$ is sub-additive, i.e. $\mathcal{H}_\mu^{q,t}(\displaystyle\bigcup_nA_n)\leq\displaystyle\sum_n\mathcal{H}_\mu^{q,t}(A_n)$.
\end{description}
The first item is obvious. Let us prove (ii). Let $E\subseteq\,F$ be nonempty subsets of $\mathbb{R}^d$. We have
$$
\mathcal{H}_\mu^{q,t}(E)=\displaystyle\sup_{A\subseteq\,E}{\overline\mathcal{H}}_\mu^{q,t}(A)
\leq\displaystyle\sup_{A\subseteq\,F}{\overline\mathcal{H}}_\mu^{q,t}(A)=\mathcal{H}_\mu^{q,t}(F).
$$
We next prove (iii). If the right hand term is infinite, the inequality is obvious. So, assume that it is finite. Let $(E_n)_n$ be a countable family of subsets $E_i\subseteq\mathbb{R}^d$ for which $\displaystyle\sum_n\mathcal{H}_\mu^{q,t}(E_n)<\infty$. Let also $\epsilon,\delta>0$ and $(B(x_{ni},r_{ni}))_i$ a centered $\epsilon$-covering of $E_n$ satisfying
$$
\displaystyle\sum_i\mu(B(x_{ni},r_{ni}))^q(2r_{ni})^t\leq{\overline\mathcal{H}}_{\mu,\delta}^{q,t}(E_n)+\frac{\delta}{2^n}.
$$
The whole set $(B(x_{ni},r_{ni}))_{n,i}$ is a centered $\epsilon$-covering of the whole union $\displaystyle\bigcup_nE_n$. As a consequence,
$$
\matrix{\displaystyle{\overline\mathcal{H}}_{\mu,\epsilon}^{q,t}(\bigcup_nE_n)
&\leq&\displaystyle\sum_n\sum_i(\mu(B(x_{ni},r_{ni})))^q(2r_{ni})^t\hfill\cr\medskip
&\leq&\displaystyle\sum_n\bigl({\overline\mathcal{H}}_{\mu,\epsilon}^{q,t}(E_n)+\frac{\delta}{2^n}\bigr)\hfill\cr\medskip
&\leq&\displaystyle\sum_n\bigl({\overline\mathcal{H}}_\mu^{q,t}(E_n)+\frac{\delta}{2^n}\bigr)\hfill\cr\medskip
&\leq&\displaystyle\sum_n\mathcal{H}_\mu^{q,t}(E_n)+\delta.\hfill}
$$
Having $\epsilon$ and $\delta$ going towards 0, we obtain
$$
{\overline\mathcal{H}}_\mu^{q,t}(\displaystyle\bigcup_nE_n)\leq\displaystyle\sum_n\mathcal{H}_\mu^{q,t}(E_n).
$$
Let next a set $F$ covered with the countable set $(A_n)_n$. That is $F\subseteq\displaystyle\bigcup_nA_n$. We have
$$
\matrix{{\overline\mathcal{H}}_\mu^{q,t}(F)
&=&{\overline\mathcal{H}}_\mu^{q,t}\biggl(\displaystyle\bigcup_n(A_n\cap\,F)\biggr)\hfill\cr\medskip
&\leq&\displaystyle\sum_n\mathcal{H}_\mu^{q,t}(A_n\cap\,F)\hfill\cr\medskip
&\leq&\displaystyle\sum_n\mathcal{H}_\mu^{q,t}(A_n).\hfill}
$$
Taking the sup on $F$, we obtain
$$
\mathcal{H}_\mu^{q,t}(F)\leq\displaystyle\sum_n\mathcal{H}_\mu^{q,t}(A_n).
$$
We now prove that $\mathcal{H}_\mu^{q,t}$ is metric. Let $A,B$ subsets of $\mathbb{R}^d$ where the distance $d(A,B)$ is defined by $d(A,B)=inf\{|x-y|;\;\;x\in\,A\,\;\;y\in\,B\}>0$ and $\mathcal{H}_\mu^{q,t}(A\cup\,B)<\infty$. Let next $0<\delta<d(A,B)$, $\varepsilon>0$, $F_1\subseteq\,A$,
$F_2\subseteq\,B$ and $(B(x_i,r_i))_i$ a centered $\delta$-covering of the set $F_1\cup\,F_2$ and such that
$$
{\overline\mathcal{H}}_{\mu,\delta}^{q,t}(F_1\cup\,F_2)\leq\displaystyle\sum_i(\mu(B(x_i,r_i)))^q(2r_i)^t\leq{\overline\mathcal{H}}_{\mu,\delta}^{q,t}(F_1\cup\,F_2)+\varepsilon.
$$
This is always possible from the definition of ${\overline\mathcal{H}}_{\mu,\delta}^{q,t}(F_1\cup\,F_2)$.
Denote next the index sets
$$
I=\{\,i;\;\;B(x_i,r_i)\cap\,F_1\not=\emptyset\,\}\quad\hbox{and}\quad\,J=\{\,i;\;\;B(x_i,r_i)\cap\,F_2\not=\emptyset\,\}.
$$
Hence, the countable sets $(B(x_i,r_i))_{i\in\,I}$ and $(B(x_i,r_i))_{i\in\,J}$ are centered $\delta$-coverings of $F_1$ and $F_2$ respectively. Consequently,
$$
\matrix{{\overline\mathcal{H}}_{\mu,\delta}^{q,t}(F_1)+{\overline\mathcal{H}}_{\mu,\delta}^{q,t}(F_2)
&\leq&\displaystyle\sum_{i\in\,I}(\mu(B(x_i,r_i)))^q(2r_i)^t+\displaystyle\sum_{i\in\,J}\mu(B(x_i,r_i))^q(2r_i)^t\hfill\cr\medskip
&=&\displaystyle\sum_i(\mu(B(x_i,r_i)))^q(2r_i)^t\hfill\cr\medskip
&\leq&{\overline\mathcal{H}}_{\mu,\delta}^{q,t}(F_1\cup\,F_2)+\varepsilon.\hfill}
$$
As a result,
$$
{\overline\mathcal{H}}_\mu^{q,t}(F_1)+{\overline\mathcal{H}}_\mu^{q,t}(F_2)\leq{\overline\mathcal{H}}_\mu^{q,t}(F_1\cup\,F_2)
+\varepsilon\leq\mathcal{H}_\mu^{q,t}(A\cup\,B)+\varepsilon.
$$
When $\varepsilon\downarrow0$ and taking the sup on the sets $F_1\subseteq\,A$ and $F_2\subseteq\,B$, we obtain
$$
\mathcal{H}_\mu^{q,t}(A\cup\,B)\geq\mathcal{H}_\mu^{q,t}(A)+\mathcal{H}_\mu^{q,t}(B).
$$
The inequality
$$
\mathcal{H}_\mu^{q,t}(A\cup\,B)\leq\mathcal{H}_\mu^{q,t}(A)+\mathcal{H}_\mu^{q,t}(B).
$$
follows from the sub-additivity property of the measure $\mathcal{H}_\mu^{q,t}$.
\begin{definition}
The restriction of $\mathcal{H}_\mu^{q,t}$ on Borel sets is called the mixed generalized Hausdorff measure on $\mathbb{R}^d$.
\end{definition}
Now, we define the mixed generalized multifractal packing measure. We use already the same notations as previously. Let
$$
{\overline\mathcal{P}}_{\mu,\epsilon}^{q,t}(E)=\displaystyle\sup\{\,\sum_i(\mu(B(x_i,r_i)))^q(2r_i)^t\,\}
$$
where the sup is taken over the set of all centered $\epsilon$-packings of $E$. For the empty set, we set as usual ${\overline\mathcal{P}}_{\mu,\epsilon}^{q,t}(\emptyset)=0$. Next, we consider the limit as $\epsilon\downarrow0$,
$$
{\overline\mathcal{P}}_\mu^{q,t}(E)=\displaystyle\lim_{\epsilon\downarrow0}{\overline\mathcal{P}}_{\mu,\epsilon}^{q,t}(E)
=\inf_{\delta>0}{\overline\mathcal{P}}_{\mu,\epsilon}^{q,t}(E)
$$
and finally,
$$
\mathcal{P}_\mu^{q,t}(E)=\displaystyle\inf_{E\subseteq\,\cup_iE_i}\sum_i{\overline\mathcal{P}}_\mu^{q,t}(E_i).
$$
\begin{lemma}\label{anouarlemme2}
$\mathcal{P}_\mu^{q,t}$ is an outer metric measure on $\mathbb{R}^d$.
\end{lemma}
The proof of this lemma uses the following result.
\begin{equation}\label{anouarclaim1}
{\overline\mathcal{P}}_\mu^{q,t}(A\cup\,B)={\overline\mathcal{P}}_\mu^{q,t}(A)+{\overline\mathcal{P}}_\mu^{q,t}(B),\;\;\hbox{whenever}\;d(A,B)>0.
\end{equation}
Indeed, let $0<\epsilon<\displaystyle\frac{1}{2}d(A,B)$ and $(B(x_i,r_i))_i$ be a centered $\epsilon$-packing of the union $A\cup\,B$. It  \\
can be divided into two parts $I$ and $J$,
$$
(B(x_i,r_i))_i=\Bigl(B(x_i,r_i)\Bigr)_{i\in\,I}\bigcup\Bigl(B(x_i,r_i)\Bigr)_{i\in\,J}
$$
where
$$
\forall\,i\in\,I,\;\;B(x_i,r_i)\cap\,B=\emptyset\quad\hbox{and}\quad\forall\,i\in\,J,\;\;B(x_i,r_i)\cap\,A=\emptyset.
$$
Therefore, $(B(x_i,r_i))_{i\in\,I}$ is a centered $\epsilon$-packing of $A$ and $(B(x_i,r_i))_{i\in\,J}$ is a centered $\epsilon$-packing of the union $B$. Hence,
$$
\displaystyle\sum_i(\mu(B(x_i,r_i)))^q(2r_i)^t=
\underbrace{\displaystyle\sum_{i\in\,I}(\mu(B(x_i,r_i)))^q(2r_i)^t}_{\leq{\overline\mathcal{P}}_{\mu,\epsilon}^{q,t}(A)}
+\underbrace{\displaystyle\sum_{i\in\,I}(\mu(B(x_i,r_i)))^q(2r_i)^t}_{\leq{\overline\mathcal{P}}_{\mu,\epsilon}^{q,t}(B)}
$$
Consequently,
$$
{\overline\mathcal{P}}_{\mu,\epsilon}^{q,t}(A\cup\,B)\leq{\overline\mathcal{P}}_{\mu,\epsilon}^{q,t}(A)+{\overline\mathcal{P}}_{\mu,\epsilon}^{q,t}(B)
$$
and thus the limit for $\epsilon\downarrow0$ gives
$$
{\overline\mathcal{P}}_{\mu}^{q,t}(A\cup\,B)\leq{\overline\mathcal{P}}_{\mu}^{q,t}(A)+{\overline\mathcal{P}}_{\mu}^{q,t}(B).
$$
The converse is more easier and it states that ${\overline\mathcal{P}}_{\mu,\epsilon}^{q,t}$ and next ${\overline\mathcal{P}}_{\mu}^{q,t}$ are sub-additive. Let $(B(x_i,r_i))_i$ be a centered $\epsilon$-packing of $A$ and $(B(y_i,r_i))_i$ be a centered $\epsilon$-packing of $B$. The union $\Bigl(B(x_i,r_i)\Bigr)_i\bigcup\Bigl(B(y_i,r_i)\Bigr)_i$ is a centered $\epsilon$-packing of $A\cup\,B$. So that
$$
{\overline\mathcal{P}}_{\mu,\epsilon}^{q,t}(A\cup\,B)\geq\displaystyle\sum_i(\mu(B(x_i,r_i)))^q(2r_i)^t+\displaystyle\sum_i(\mu(B(y_i,r_i)))^q(2r_i)^t.
$$
Taking the sup on $(B(x_i,r_i))_i$ as a centered $\epsilon$-packing of $A$ and next the sup on $(B(y_i,r_i))_i$ as a centered $\epsilon$-packing of $B$, we obtain
$$
{\overline\mathcal{P}}_{\mu,\epsilon}^{q,t}(A\cup\,B)\geq{\overline\mathcal{P}}_{\mu,\epsilon}^{q,t}(A)+{\overline\mathcal{P}}_{\mu,\epsilon}^{q,t}(B)
$$
and thus the limit for $\epsilon\downarrow0$ gives
$$
{\overline\mathcal{P}}_{\mu}^{q,t}(A\cup\,B)\geq{\overline\mathcal{P}}_{\mu}^{q,t}(A)+{\overline\mathcal{P}}_{\mu}^{q,t}(B).
$$
{\bf Proof of Lemma \ref{anouarlemme2}.} We shall prove as previously
\begin{description}
\item [i.] $\mathcal{P}_\mu^{q,t}(\emptyset)=0$.
\item[ii.] $\mathcal{P}_\mu^{q,t}$ is monotone, i.e. $\mathcal{P}_\mu^{q,t}(E)\leq\mathcal{P}_\mu^{q,t}(F)$, whenever $E\subseteq\,F\subseteq\mathbb{R}^d$.
\item[iii.] $\mathcal{P}_\mu^{q,t}$ is sub-additive, i.e. $\mathcal{P}_\mu^{q,t}(\displaystyle\bigcup_nA_n)\leq\displaystyle\sum_n\mathcal{P}_\mu^{q,t}(A_n)$.
\end{description}
The first item is immediate from the definition of $\mathcal{P}_\mu^{q,t}(\emptyset)=0$. Let $E\subseteq\,F$ be subsets of $\mathbb{R}^d$. We have
$$
\mathcal{P}_\mu^{q,t}(E)=\displaystyle\inf_{E\subseteq\displaystyle\bigcup_iE_i}\displaystyle\sum_i{\overline\mathcal{P}}_\mu^{q,t}(E_i)
\leq\displaystyle\inf_{F\subseteq\displaystyle\bigcup_iE_i}\displaystyle\sum_i{\overline\mathcal{P}}_\mu^{q,t}(E_i)
=\mathcal{P}_\mu^{q,t}(F).
$$
So is the item {\bf ii}. Let next $(A_n)_n$ a countable set of subsets of $\mathbb{R}^d$, $\varepsilon>0$ and for each $n$, $(E_{ni})_i$ be a covering of
$A_n$ such that
$$
\displaystyle\sum_i{\overline\mathcal{P}}_\mu^{q,t}(E_{ni})\leq\mathcal{P}_\mu^{q,t}(A_n)+\frac{\varepsilon}{2^n}.
$$
It follows for all $\varepsilon>0$ that
$$
\mathcal{P}_\mu^{q,t}(\displaystyle\bigcup_nA_n)\leq\displaystyle\sum_n
\displaystyle\sum_i{\overline\mathcal{P}}_\mu^{q,t}(E_{ni})
\leq\displaystyle\sum_n\mathcal{P}_\mu^{q,t}(A_n)+\varepsilon.
$$
Hence,
$$
\mathcal{P}_\mu^{q,t}(\displaystyle\bigcup_nA_n)\leq\displaystyle\sum_n\mathcal{P}_\mu^{q,t}(A_n).
$$
So is the item {\bf iii}. We now prove that $\mathcal{P}_\mu^{q,t}$ is metric. Let $A,B$ subsets of $\mathbb{R}^d$ be such that $d(A,B)>0$. We shall prove that $$
\mathcal{P}_\mu^{q,t}(A\cup\,B)=\mathcal{P}_\mu^{q,t}(A)+\mathcal{P}_\mu^{q,t}(B).
$$
Since $\mathcal{P}_\mu^{q,t}$ is an outer measure, it suffices to show that
$$
\mathcal{P}_\mu^{q,t}(A\cup\,B)\geq\mathcal{P}_\mu^{q,t}(A)+\mathcal{P}_\mu^{q,t}(B).
$$
Of course, if the left hand term is infinite, the inequality is obvious. So, suppose that it is finite. For $\varepsilon>0$, there exists a covering $(E_i)_i$ of the union set $A\cup\,B$ such that
$$
\displaystyle\sum_i{\overline\mathcal{P}}_\mu^{q,t}(E_i)\leq\mathcal{P}_\mu^{q,t}(A\cup\,B)+\varepsilon.
$$
By denoting $F_i=A\cap\,E_i$ and $H_i=B\cap\,E_i$, we get countable coverings $(F_i)_i$ of $A$ and $(H_i)_i$ for $B$ respectively. Furthermore, $F_i\cap\,H_j=\emptyset$ pour all $i$ and $j$. Consequently,
$$
\mathcal{P}_\mu^{q,t}(A)+\mathcal{P}_\mu^{q,t}(B)\leq\displaystyle\sum_i({\overline\mathcal{P}}_\mu^{q,t}(F_i)+{\overline\mathcal{P}}_\mu^{q,t}(H_i)).
$$
Since $d(A,B)>0$, $F_i\subset\,A$ and $H_i\subset\,B$, it follows that $d(F_i,H_j)>0$ for all $i$ and $j$. Hence, claim \ref{anouarclaim1} affirms that
$$
{\overline\mathcal{P}}_\mu^{q,t}(E_i)={\overline\mathcal{P}}_\mu^{q,t}(F_i\cup\,H_i)={\overline\mathcal{P}}_\mu^{q,t}(F_i)+{\overline\mathcal{P}}_\mu^{q,t}(H_i).
$$
Hence,
$$
\mathcal{P}_\mu^{q,t}(A)+\mathcal{P}_\mu^{q,t}(B)\leq\displaystyle\sum_i{\overline\mathcal{P}}_\mu^{q,t}(E_i)\leq\mathcal{P}_\mu^{q,t}(A\cup\,B)+\varepsilon
$$
and the result is obtained by having $\varepsilon\downarrow0$.
\begin{definition}
The restriction of $\mathcal{P}_\mu^{q,t}$ on Borel sets is called the mixed generalized Hausdorff measure on $\mathbb{R}^d$.
\end{definition}
It holds as for the case of the multifractal analysis of a single measure that the measures $\mathcal{H}_{\mu}^{q,t}$, $\mathcal{P}_\mu^{q,t}$ and the pre-measure ${\overline\mathcal{P}}_\mu^{q,t}$ assign a dimension to every set $E\subseteq\mathbb{R}^d$.
\begin{proposition}\label{anouarproposition1}
Given a subset $E\subseteq\mathbb{R}^d$,
\begin{enumerate}
\item There exists a unique number $dim_\mu^q(E)\in[-\infty,+\infty]$ such that
$$
\mathcal{H}_\mu^{q,t}(E)=\left\{\matrix{+\infty&\hbox{for}&t<dim_\mu^q(E)\hfill\cr\medskip
0&\hbox{si}&t>dim_\mu^q(E)\hfill}\right.
$$
\item There exists a unique number $Dim_\mu^q(E)\in[-\infty,+\infty]$ such that
$$
\mathcal{P}_\mu^{q,t}(E)=\left\{\matrix{+\infty&\hbox{for}&t<Dim_\mu^q(E)\hfill\cr\medskip
0&\hbox{for}&t>Dim_\mu^q(E)\hfill}\right.
$$
\item There exists a unique number $\Delta_\mu^q(E)\in[-\infty,+\infty]$ such that
$$
{\overline\mathcal{P}}_\mu^{q,t}(E)=\left\{\matrix{+\infty&\hbox{for}&t<\Delta_\mu^q(E)\hfill\cr\medskip
0&\hbox{for}&t>\Delta_\mu^q(E)\hfill}\right.
$$
\end{enumerate}
\end{proposition}
\begin{definition}\label{mixeddimensions}
The quantities $dim_\mu^q(E)$, $Dim_\mu^q(E)$ and $\Delta_\mu^q(E)$ defines the so-called mixed multifractal generalizations of the Hausdorff dimension, the packing dimension and the logarithmic index of the set $E$.
\end{definition}
Remark that if we denote $1_i=(0,0,...,q_i,0,...,0)$ the vector with zero coordinates except the ith one which equals 1, we obtain the multifractal generalizations of the Hausdorff dimension, the packing dimension and the logarithmic index of the set $E$ for the single measure $\mu_i$,
$$
dim_\mu^{1_i}(E)=dim_{\mu_i}^{q_i}(E),\quad\,Dim_\mu^{1_i}(E)=Dim_{\mu_i}^{q_i}(E)\quad\hbox{and}\quad\Delta_\mu^{1_i}(E)=\Delta_{\mu_i}^{q_i}(E).
$$
Similarly, for the null vector of $\mathbb{R}^k$, we obtain
$$
dim_\mu^{0}(E)=dim(E),\quad\,Dim_\mu^{0}(E)=Dim(E)\quad\hbox{and}\quad\Delta_\mu^{0}(E)=\Delta(E).
$$
{\bf Proof of Proposition \ref{anouarproposition1}.}\\
1. We claim that $\forall\,t\in\mathbb{R}$ such that $\mathcal{H}_\mu^{q,t}(E)<\infty$ it holds that $\mathcal{H}_\mu^{q,t'}(E)=0$ for any $t'>t$. Indeed, let $\epsilon>0$, $F\subseteq\,E$ and $(B(x_i,r_i))_i$ be a centered $\epsilon$-covering of $F$. We have
$$
{\overline\mathcal{H}}_{\mu,\epsilon}^{q,t'}(F)\leq\displaystyle\sum_i(\mu(B(x_i,r_i)))^q(2r_i)^{t'}
\leq\delta^{t'-t}\displaystyle\sum_i(\mu(B(x_i,r_i)))^q(2r_i)^t.
$$
Consequently,
$$
{\overline{H}}_{\mu,\epsilon}^{q,t'}(F)\leq\epsilon^{t'-t}{\overline{H}}_{\mu,\epsilon}^{q,t'}(F).
$$
Hence,
$$
{\overline\mathcal{H}}_\mu^{q,t'}(F)=0,\quad\forall\,F\subseteq\,E.
$$
As a result, $\mathcal{H}_\mu^{q,t'}(E)=0$. We then set
$$
dim_{\mu}^q(E)=\inf\{\,t\in\mathbb{R};\;\;\mathcal{H}_\mu^{q,t'}(E)=0\,\}.
$$
One can proceed otherwise by claiming that $\forall\,t\in\mathbb{R}$ such that $\mathcal{H}_\mu^{q,t}(E)>0$ it holds that $\mathcal{H}_\mu^{q,t'}(E)=+\infty$ for any $t'<t$. Indeed, proceeding as previously, we obtain for $\epsilon>0$,
$$
\epsilon^{t'-t}{\overline{H}}_{\mu,\epsilon}^{q,t}(F)\leq{\overline{H}}_{\mu,\epsilon}^{q,t'}(F).
$$
Hence,
$$
{\overline\mathcal{H}}_\mu^{q,t'}(F)=+\infty,\quad\forall\,F\subseteq\,E.
$$
As a result, $\mathcal{H}_\mu^{q,t'}(E)=+\infty$. We then set
$$
dim_{\mu}^q(E)=\sup\{\,t\in\mathbb{R};\;\;\mathcal{H}_\mu^{q,t'}(E)=+\infty\,\}.
$$
2. Similarly to the previous case, let $t\in\mathbb{R}$ be such that $\mathcal{P}_\mu^{q,t}(E)<\infty$. There exists $(E_i)_i$ subsets of $\mathbb{R}^d$ satisfying
$$
E\subseteq\displaystyle\bigcup_i{E_i},\quad\hbox{and}\quad{\overline{P}}_\mu^{q,t}(E_i)<\infty,\;\;\hbox{for any}\,i.
$$
Let next $t'>t$, $\epsilon>0$ and $(B(x_{ni},r_{ni}))_n$ be a centered $\epsilon$-packing of the set $E_i$. Then
$$
\displaystyle\sum_n(\mu(B(x_{ni},r{ni})))^q(2r_{ni})^{t'}\leq\epsilon^{t'-t}\displaystyle\sum_n(\mu(B(x_{ni},r{ni})))^q(2r_{ni})^t
$$
which implies that
\begin{equation}\label{Pmuqtcoupure}
{\overline{P}}_{\mu,\epsilon}^{q,t'}(E_i)\leq\epsilon^{t'-t}\overline{P}_{\mu,\epsilon}^{q,t}(E_i).
\end{equation}
Hence, ${\overline{P}}_\mu^{q,t'}(E_i)=0$ for all $i$ and consequently $\mathcal{P}_\mu^{q,t'}(E)=0$. we set as previously
$$
Dim_{\mu}^q(E)=\inf\{\,t\in\mathbb{R};\;\;\mathcal{P}_\mu^{q,t}(E)=0\,\}.
$$
3. It follows from equation \ref{Pmuqtcoupure} that for any $t\in\mathbb{R}$ such that ${\overline\mathcal{P}}_\mu^{q,t}(E)<\infty$, we have ${\overline\mathcal{P}}_\mu^{q,t'}(E)=0$ for any $t'>t$. We then set
$$
\Delta_{\mu}^q=\inf\{\,t\in\mathbb{R};\;\;\mathcal{P}_\mu^{q,t}(E)=0\,\}.
$$
Next, we aim to study the characteristics of the mixed multifractal generalizations of dimensions. To do this we will adapt the following notations. For $q=(q_1,...,q_k)\in\mathbb{R}^k$,
$$
b_{\mu,E}(q)=dim_{\mu}^{q}(E),\,\,B_{\mu,E}(q)=Dim_{\mu}^{q}(E)\,\,\hbox{and}\,\,\Lambda_{\mu,E}(q)=\Delta_{\mu}^{q}(E).
$$
When $E=support(\mu)$ is the support of the measure $\mu$, we will omit the indexation with $E$ and denote simply
$$
b_{\mu}(q),\,\,B_{\mu}(q)\,\,\hbox{and}\,\,\Lambda_{\mu}(q).
$$
The following propositions resumes the characteristics of these functions and extends the results of L. Olsen \cite{olsen1} for our case.
\begin{proposition}\label{anouarproposition2}
\begin{description}
\item[a.] $0\leq\,b_{\mu,E}(q)$, whenever $q_i\leq1$ and $\mu_i(E)>0$ for all $i=1,2,...,k$.
\item[b.] $\Lambda_{\mu,E}(q)\leq0$, whenever $q_i\geq1$ for all $i=1,2,...,k$.
\item[c.] $b_{\mu,.}(q)$ and $B_{\mu,.}(q)$ are non decreasing with respect to the inclusion property in $\mathbb{R}^d$.
\item[c.] $b_{\mu,.}(q)$ and $B_{\mu,.}(q)$ are $\sigma$-stable.
\end{description}
\end{proposition}
{\bf Proof.}
{\bf a.} ...\\
{\bf b.} For $q\leq1$, it holds that
$$
{\overline\mathcal{P}}_{\mu,\epsilon}^{q,t}(E)\leq\epsilon^t,\;\;\forall\,t>0.
$$
Hence,
$$
{\overline\mathcal{P}}_{\mu}^{q,t}(E)=0,\;\;\forall\,t>0
$$
which means that
$$
\Lambda_{\mu,E}(q)\leq\,t,\;\;\forall\,t>0\quad\Longleftrightarrow\quad\Lambda_{\mu,E}(q)\leq0.
$$
{\bf c.} Let $E\subseteq\,F$ be subsets of $\mathbb{R}^d$. We have
$$
\mathcal{H}_{\mu}^{q,t}(E)=\displaystyle\sup_{A\subseteq\,E}{\overline\mathcal{H}}_{\mu}^{q,t}(A)
\leq\displaystyle\sup_{A\subseteq\,F}{\overline\mathcal{H}}_{\mu}^{q,t}(A)=\mathcal{H}_{\mu}^{q,t}(F).
$$
So for the monotony of $b_{\mu,.}(q)$. Next, since $\mathcal{P}_\mu^{q,t}$ is an outer measure,
$$
\mathcal{P}_\mu^{q,t}(E)\leq\mathcal{P}_\mu^{q,t}(F)=0,\quad\forall\,t>B_{\mu,F}(q).
$$
Consequently,
$$
\mathcal{P}_\mu^{q,t}(E)=0,\quad\forall\,t>B_{\mu,F}(q).
$$
Therefore,
$$
B_{\mu,E}(q)\leq t,,\quad\forall\,t>B_{\mu,F}(q).
$$
So that
$$
B_{\mu,E}(q)\leq\,B_{\mu,F}(q).
$$
{\bf d.} Let $(A_n)_n$ be a countable set of subsets $A_n\subseteq\mathbb{R}^d$ and denote $A=\displaystyle\bigcup_nA_n$. It holds from the monotony of $b_{\mu,.}(q)$ that
$$
b_{\mu,A_n}(q)\leq\,b_{\mu,A}(q),\quad\forall\,n.
$$
Hence,
$$
\displaystyle\sup_nb_{\mu,A_n}(q)\leq\,b_{\mu,A}(q).
$$
Next, for any $t>\displaystyle\sup_nb_{\mu,A_n}(q)$, there holds that
$$
\mathcal{H}_\mu^{q,t}(A_n)=0,\quad\forall\,n.
$$
Consequently, from the sub-additivity property of $\mathcal{H}_\mu^{q,t}$, it holds that
$$
\mathcal{H}_\mu^{q,t}(\displaystyle\bigcup_nA_n)=0,\quad\forall\,t>\displaystyle\sup_nb_{\mu,A_n}(q).
$$
Which means that
$$
b_{\mu,A}(q)\leq\,t,\quad\forall\,t>\displaystyle\sup_nb_{\mu,A_n}(q).
$$
Hence,
$$
b_{\mu,A}(q)\leq\displaystyle\sup_nb_{\mu,A_n}(q).
$$
We shall now prove the $\sigma$-stability of $B_{\mu,.}(q)$. Consider as previously a countable set $(A_n)_n$ of subsets of $\mathbb{R}^d$. The following inequality is immediate.
$$
\displaystyle\sup_nB_{\mu,A_n}(q)\leq\,B_{\mu,A}(q).
$$
Next, for any $t>\displaystyle\sup_nB_{\mu,A_n}(q)$, we have
$$
\mathcal{P}_\mu^{q,t}(A)\leq\displaystyle\sum_n\mathcal{P}_\mu^{q,t}(A_n)=0.
$$
So that
$$
\mathcal{P}_\mu^{q,t}(A)=0,\quad\forall\,t>\displaystyle\sup_nB_{\mu,A_n}(q).
$$
Consequently,
$$
B_{\mu,A}(q)\leq\,t,\quad\forall\,t>\displaystyle\sup_nB_{\mu,A_n}(q).
$$
Which means that
$$
B_{\mu,A}(q)\leq\displaystyle\sup_nB_{\mu,A_n}(q).
$$
Next, we continue to study the characteristics of the mixed generalized multifractal dimensions. The following result is obtained.
\begin{proposition}\label{anouarproposition3}
\begin{description}
\item[a.] The functions $q\longmapsto\,B_{\mu}(q)$ and $q\longmapsto\Lambda_{\mu}(q)$ are convex.
\item[b.] For $i=1,2,...,k$, the functions $q_i\longmapsto\,b_{\mu}(q)$, $q_i\longmapsto\,B_{\mu}(q)$ and $q_i\longmapsto\Lambda_{\mu}(q)$, ($\widehat{q_i}=(q_1,\dots,q_{i-1},q_{i+1},\dots,q_k)$ fixed), are non increasing.
\end{description}
\end{proposition}
{\bf Proof.}
{\bf a.} We start by proving that $\Lambda_{\mu,E}$ is convex. Let $p,\,q\in\mathbb{R}^k$, $\alpha\in]0,1[$, $s>\Lambda_{\mu,E}(p)$ and $t>\Lambda_{\mu,E}(q)$. Consider next a centered $\epsilon$-packing $(B_i=B(x_i,r_i))_i$ of $E$. Applying H\"older's inequality, it holds that
$$
\displaystyle\sum_i(\mu(B_i))^{\alpha\,q+(1-\alpha)p}(2r_i)^{\alpha\,t+(1-\alpha)s}\leq\biggl(\displaystyle\sum_i(\mu(B_i))^q(2r_i)^t\biggr)^\alpha
\biggl(\displaystyle\sum_i(\mu(B_i))^p(2r_i)^s\biggr)^{1-\alpha}.
$$
Hence,
$$
{\overline\mathcal{P}}_{\mu,\epsilon}^{\alpha\,q+(1-\alpha)p,\alpha\,t+(1-\alpha)s}(E)\leq\biggl({\overline\mathcal{P}}_{\mu,\epsilon}^{q,t}(E)\biggr)^\alpha
\biggl({\overline\mathcal{P}}_{\mu,\epsilon}^{p,s}(E)\biggr)^{1-\alpha}.
$$
The limit on $\epsilon\downarrow0$ gives
$$
{\overline\mathcal{P}}_\mu^{\alpha q+(1-\alpha)p,\alpha\,t+(1-\alpha)s}(E)\leq\biggl({\overline\mathcal{P}}_\mu^{q,t}(E)\biggr)^\alpha
\biggl({\overline\mathcal{P}}_\mu^{p,s}(E)\biggr)^{1-\alpha}.
$$
Consequently,
$$
{\overline\mathcal{P}}_\mu^{\alpha q+(1-\alpha)p,\alpha\,t+(1-\alpha)s}(E)=0\,\;\;\forall\,s>\Lambda_{\mu,E}(p)\;\;\hbox{and}\;\;t>\Lambda_{\mu,E}(q).
$$
It results that
$$
\Lambda_{\mu,E}(\alpha q+(1-\alpha)p)\leq\alpha\Lambda_{\mu,E}(q)+(1-\alpha)\Lambda_{\mu,E}(p).
$$
We now prove the convexity of $B_{\mu,E}$. We set in this case $t=B_{\mu,E}(q)$ and $s=B_{\mu,E}(p)$. We have
$$
\mathcal{P}_\mu^{q,t+\varepsilon}(E)=\mathcal{P}_\mu^{p,s+\varepsilon}(E)=0.
$$
Therefore, there exists $(H_i)_i$ and $(K_i)_i$ coverings of the set $E$ for which
$$
\displaystyle\sum_i{\overline\mathcal{P}}_\mu^{q,t+\varepsilon}(H_i)\leq1
\qquad\hbox{et}\qquad
\displaystyle\sum_i{\overline\mathcal{P}}_\mu^{p,s+\varepsilon}(K_i)\leq1.
$$
Denote for $n\in\mathbb{N}$, $E_n=\displaystyle\bigcup_{1\leq\,i,j\leq\,n}(H_i\cap K_j)$. Thus, $(E_n)_n$ is a covering of $E$. So that,
$$
\matrix{& &\mathcal{P}_\mu^{\alpha\,q+(1-\alpha)p,\alpha\,t+(1-\alpha)s+\varepsilon}(E_n)\hfill\cr\medskip
&\leq&\displaystyle\sum_{i,j=1}^n\mathcal{P}_\mu^{\alpha\,q+(1-\alpha)p,\alpha t+(1-\alpha)s+\varepsilon}(H_i\cap\,K_j)\hfill\cr\medskip
&\leq&\displaystyle\sum_{i,j=1}^n{\overline\mathcal{P}}_\mu^{\alpha\,q+(1-\alpha)p,\alpha t+(1-\alpha)s+\varepsilon}(H_i\cap\,K_j)\hfill\cr\medskip
&\leq&\biggl(\displaystyle\sum_{i,j=1}^n{\overline\mathcal{P}}_\mu^{q,t+\varepsilon}(H_i\cap\,K_j)\biggr)^\alpha
\biggl(\displaystyle\sum_{i,j=1}^n{\overline\mathcal{P}}_\mu^{p,s+\varepsilon}(H_i\cap\,K_j)\biggr)^{1-\alpha}\hfill\cr\medskip
&\leq&n^\alpha\,n^{1-\alpha}=n<\infty.\hfill}
$$
Consequently,
$$
B_{\mu,E_n}(\alpha\,q+(1-\alpha)p)\leq\alpha\,t+(1-\alpha)s+\varepsilon,\quad\forall\,\varepsilon>0.
$$
Hence,
$$
B_{\mu,E}(\alpha\,q+(1-\alpha)p)\leq\alpha\,B_{\mu,E}(q)+(1-\alpha)B_{\mu,E}(p).
$$
{\bf b.} For $i=1,2,\dots,k$, let $\widehat{q_i}$ fixed and $p_i\leq\,q_i$ reel numbers. Denote next $q=(q_1,\dots,q_{i-1},q_i,q_{i+1},\dots,q_k)$ and $p=(q_1,\dots,q_{i-1},p_i,q_{i+1},\dots,q_k)$. Let finally $A\subseteq\,E$. For a centered $\epsilon$-covering $(B(x_i,r_i))_i$ of $A$, we have immediately
$$
\mu(B(x_i,r_i))^{q}(2r_i)^t\leq\mu(B(x_i,r_i))^{p}(2r_i)^t,\;\;\forall\,t\in\mathbb{R}.
$$
Hence,
$$
{\overline{H}}_{\mu,\epsilon}^{q,t}(A)\leq{\overline{H}}_{\mu,\epsilon}^{p,t}(A).
$$
When $\epsilon\downarrow0$, we obtain
$$
{\overline{H}}_\mu^{q,t}(A)\leq{\overline{H}}_\mu^{p,t}(A).
$$
Therefore,
$$
\mathcal{H}_\mu^{q,t}(E)=\displaystyle\sup_{A\subseteq\,E}{\overline\mathcal{H}}_\mu^{q,t}(A)
\leq\displaystyle\sup_{A\subseteq\,E}{\overline\mathcal{H}}_\mu^{p,t}(A)={\mathcal{H}}_\mu^{p,t}(E).
$$
This induces the fact that
$$
{\mathcal{H}}_\mu^{q,t}(E)=0,\quad\forall\,t>b_{\mu,E}(p).
$$
Consequently
$$
b_{\mu,E}(q)<t,\quad\forall\,t>b_{\mu,E}(p).
$$
Hence,
$$
b_{\mu,E}(q)\leq\,b_{\mu,E}(p).
$$
We shall now prove the monotony $\Lambda_{\mu,E}$. With the same notations as above and using a centered $\epsilon$-packing of the set $E$, we obtain
$$
{\overline{P}}_\mu^{q,t}(E)\leq{\overline{P}}_\mu^{p,t}(E).
$$
As a consequence,
$$
{\overline{P}}_\mu^{q,t}(E)=0,\quad\forall\,t>\Lambda_{\mu,E}(p).
$$
Therefore,
$$
\Lambda_{\mu,E}(q)<t,\quad\forall\,t>\Lambda_{\mu,E}(p).
$$
Hence,
$$
\Lambda_{\mu,E}(q)\leq\Lambda_{\mu,E}(p).
$$
We now prove that $B_{\mu,E}$ is non increasing. For $i=1,2,\dots,k$, let $\widehat{q_i}$ fixed and $p_i\leq\,q_i$ reel numbers. Denote next $q=(q_1,\dots,q_{i-1},q_i,q_{i+1},\dots,q_k)$ and $p=(q_1,\dots,q_{i-1},p_i,q_{i+1},\dots,q_k)$. Let next $(E_i)_i$ be a covering of the set $E$. It results from the previous case that
$$
\displaystyle\sum_i{\overline\mathcal{P}}_\mu^{p,t}(E_i)\geq\displaystyle\sum_i{\overline\mathcal{P}}_\mu^{q,t}(E_i).
$$
Which means that $\mathcal{P}_\mu^{p,t}(E)\geq\mathcal{P}_\mu^{q,t}(E)$. Consequently,
$$
\mathcal{P}_\mu^{q,t}(E)=0,\quad\forall\,t>B_{\mu,E}(p)
$$
and thus
$$
B_{\mu,E}(q)<t,\quad\forall\,t>B_{\mu,E}(p).
$$
Hence,
$$
B_{\mu,E}(q)\leq\,B_{\mu,E}(p).
$$
\begin{proposition}\label{anouarproposition4}
\begin{description}
\item[a.] $0\leq\,b_\mu(q)\leq\,B_\mu(q)\leq\Lambda_\mu(q)$, whenever $q_i<1$ for all $i=1,2,...,k$.
\item[b.] $b_\mu(\P_i)=B_\mu(\P_i)=\Lambda_\mu(\P_i)=0$, where $\P_i=(0,0,...,1,0,...,0)$.
\item[c.] $b_\mu(q)\leq\,B_\mu(q)\leq\Lambda_\mu(q)\leq0$ whenever $q_i>1$ for all $i=1,2,...,k$.
\end{description}
\end{proposition}
The proof of this results reposes on the following intermediate ones.
\begin{lemma}\label{anouarlemme3}
There exists a constant $\xi\in]0,+\infty[$ satisfying for any $E\subseteq\mathbb{R}^d$,
$$
\mathcal{H}_{\mu}^{q,t}(E)\leq\xi\mathcal{P}_{\mu}^{q,t}(E)\leq\xi{\overline\mathcal{P}}_{\mu}^{q,t}(E),\qquad\forall\,q,t.
$$
More precisely, $\xi$ is the number related to the Besicovitch covering theorem.
\end{lemma}
\begin{theorem}\label{besicovitch}
There exists a constant $\xi\in\mathbb{N}$ satisfying: For any $E\in\mathbb{R}^d$ and $(r_x){x\in\,E}$ a bounded set of positive real numbers, there exists $\xi$ sets $B_1$, $B_2$, ..., $B_{\xi}$, that are finite or countable composed of balls $B(x,r_x)$, $x\in\,E$ such that
\begin{itemize}
\item $E\subseteq\displaystyle\bigcup_{1\leq\,i\leq\xi}\displaystyle\bigcup_{B\in\,B_i}B$.
\item each $B_i$ is composed of disjoint balls.
\end{itemize}
\end{theorem}
{\bf Proof of Lemma \ref{anouarlemme3}.} It suffices to prove the first inequality. The second is always true for all $\xi>0$. Let $F\subseteq\mathbb{R}^d$, $\epsilon>0$ and $\mathcal{V}=\{\,B(x,\frac{\epsilon}{2});\;\;\,x\in\,F\,\}$. Let next $\bigl((B_{ij})_j\bigr)_{1\leq\,i\leq\xi}$ be the $\xi$ sets of $\mathcal{V}$ obtained by the Besicovitch covering theorem. So that, $(B_{ij})_{i,j}$ is a centered $\epsilon$-covering of the set $F$ and for each $i$, $(B_{ij})_j$ is a centered $\epsilon$-packing of $F$. Therefore,
$$
\displaystyle{\overline\mathcal{H}}_{\mu,\epsilon}^{q,t}(F)\leq\sum_{i=1}^{\xi}\sum_j\bigl(\mu(B_{ij})\bigr)^q(2r_{ij})^t
\leq\sum_{i=1}^{\xi}{\overline\mathcal{P}}_{\mu,\epsilon}^{q,t}(F)=\xi{\overline\mathcal{P}}_{\mu,\epsilon}^{q,t}(F).
$$
Hence, ${\overline\mathcal{H}}_{\mu}^{q,t}(F)\leq\xi{\overline\mathcal{P}}_{\mu}^{q,t}(F)$. Consequently, for $E\subseteq\displaystyle\bigcup_iE_i$, we obtain
$$
\matrix{\mathcal{H}_{\mu}^{q,t}(E)=\mathcal{H}_{\mu}^{q,t}(\displaystyle\bigcup_i(E_i\cap\,E))
&\leq&\displaystyle\sum_i\mathcal{H}_{\mu}^{q,t}(E_i\cap\,E)\hfill\cr\medskip
&\leq&\displaystyle\sum_i\displaystyle\sup_{F\subseteq E_i\cap\,E}{\overline\mathcal{H}}_{\mu}^{q,t}(F)\hfill\cr\medskip
&\leq&\xi\displaystyle\sum_i\displaystyle\sup_{F\subseteq E_i\cap\,E}{\overline\mathcal{P}}_{\mu}^{q,t}(F)\hfill\cr\medskip
&\leq&\xi\displaystyle\sum_i{\overline\mathcal{P}}_{\mu}^{q,t}(E_i).\hfill}
$$
So as Lemma \ref{anouarlemme3}.\\
{\bf Proof of Proposition \ref{anouarproposition4}.} It follows from Proposition \ref{anouarproposition2}, Proposition \ref{anouarproposition3} and Lemma \ref{anouarlemme3}.
\section{Mixed multifractal generalization of Bouligand-Minkowski's dimension}
In this section, we propose to develop mixed multifractal generalization of Bouligand-Minkowski's dimension. Such a dimension is sometimes called the box-dimension or the Renyi dimension. Some mixed generalizations are already introduced in \cite{olsen3}. We will see hereafter that the mixed generalizations to be provided resemble to those in \cite{olsen3}. We will prove that in the mixed case, these dimensions remain strongly related to the mixed multifractal generalizations of the Hausdorff and packing dimensions. In the case of a single measure $\mu$, the Bouligand-Minkowski dimensions are introduced as follows. For $E\subseteq\hbox{Support}(\mu)$, $\delta>0$ and $q\in\mathbb{R}$, let
$$
\mathcal{T}_{\mu,\delta}^q(E)=\inf\left\{\displaystyle\sum_i\bigl(\mu\bigl(B(x_i,\delta)\bigr)\bigr)^q\right\}
$$
where the inf is over the set of all centered $\delta$-coverings $\bigl(B(x_i,\delta)\bigr)_i$ of the set $E$. The Bouligand-Minkowski dimensions are
$$
{\overline{L}}_{\mu}^q(E)=\displaystyle\limsup_{\delta\downarrow0}\displaystyle\frac{\log\bigl(\mathcal{T}_{\mu,\delta}^q(E)\bigr)}{-\log\delta}
$$
for the upper one and
$$
{\underline{L}}_{\mu}^q(E)=\displaystyle\liminf_{\delta\downarrow0}\displaystyle\frac{\log\bigl(\mathcal{T}_{\mu,\delta}^q(E)\bigr)}{-\log\delta}
$$
for the lower. In the case of equality, the common value is denoted ${L}_{\mu}^q(E)$ and is called the Bouligand-Minkowski dimension of the set $E$. We can equivalently define these dimensions via the $\delta$-packings as follows. For $\delta>0$ and $q\in\mathbb{R}$, we set
$$
\mathcal{S}_{\mu,\delta}^q(E)=\sup\left\{\displaystyle\sum_i\bigl(\mu\bigl(B(x_i,\delta)\bigr)\bigr)^q\right\}
$$
where the sup is taken over all the centered $\delta$-packings $\bigl(B(x_i,\delta)\bigr)_i$ of the set $E$. The upper dimension is
$$
{\overline{C}}_{\mu}^q(E)=\displaystyle\limsup_{\delta\downarrow0}\displaystyle\frac{\log\bigl(\mathcal{S}_{\mu,\delta}^q(E)\bigr)}{-\log\delta}
$$
and the lower is
$$
{\underline{C}}_{\mu}^q(E)=\displaystyle\liminf_{\delta\downarrow0}\displaystyle\frac{\log\bigl(\mathcal{S}_{\mu,\delta}^q(E)\bigr)}{-\log\delta}
$$
and similarly, when these are equal, the common value will be denoted ${C}_{\mu}^q(E)$ and it defines the dimension of $E$. We now introduce the mixed multifractal generalization of the Bouligand-Minkowski dimensions. As we have noticed, our ideas here is quite the same as the one in \cite{olsen3}. Let $\mu=(\mu_1,\,\mu_2,\dots,\mu_k)$ be a vector valued measure composed of probability measures on $\mathbb{R}^d$. Denote as previously
$$
\mu(B(x,r))\equiv\bigl(\mu_1(B(x,r)),\dots,\,\mu_k(B(x,r))\bigr)
$$
and for $q=(q_1,q_2,\dots,q_k)\in\mathbb{R}^k$,
$$
(\mu(B(x,r)))^q\equiv(\mu_1(B(x,r)))^{q_1}\dots(\mu_k(B(x,r)))^{q_k}.
$$
Next, for a nonempty subset $E\subseteq\mathbb{R}^d$ and $\delta>0$, we will use the same notations for $\mathcal{T}_{\mu,\delta}^q(E)$, ${\overline{C}}_{\mu}^q(E)$ and ${\underline{C}}_{\mu}^q(E)$ but without forgetting that we use the new product for the measure $\mu$. Similarly for
$\mathcal{S}_{\mu,\delta}^q(E)$, ${\overline{L}}_{\mu}^q(E)$ and ${\underline{L}}_{\mu}^q(E)$.
\begin{definition}\label{bouligandminkowskimixte}
For $E\subseteq\hbox{Support}(\mu)$ and $q=(q_1,q_2,\dots,q_k)\in\mathbb{R}^k$, we will call
\begin{description}
\item[a.] ${\overline{C}}_{\mu}^q(E)$ and ${\overline{L}}_{\mu}^q(E)$ the upper mixed multifractal generalizations of the Bouligand Minkowski dimension of $E$.
\item[b.] ${\underline{C}}_{\mu}^q(E)$ and ${\underline{L}}_{\mu}^q(E)$ the lower mixed multifractal generalizations of the Bouligand Minkowski dimension of $E$.
\item[c.] ${{C}}_{\mu}^q(E)$ and ${{L}}_{\mu}^q(E)$ the mixed multifractal generalizations of the Bouligand Minkowski dimension of $E$.
\end{description}
\end{definition}
\begin{remark} We stress the fact that each quantity defines in fact a mixed generalization that can be different from the other. That is, we did not mean that ${\overline{C}}_{\mu}^q(E)$ and ${\overline{L}}_{\mu}^q(E)$ are the same (equal) and similarly for the lower ones. We will prove in the contrary that as for the single case, they can be different.
\end{remark}
\begin{theorem}\label{anouartheorem1}
\begin{enumerate}
\item For all $q\in\mathbb{R}^{k}$, we have
$$
{\underline{L}}_{\mu}^q(E)\leq{\underline{C}}_{\mu}^q(E)\quad\hbox{and}\quad{\overline{L}}_{\mu}^q(E)\leq{\overline{C}}_{\mu}^q(E).
$$
\item For any $q\in\mathbb{R}^{*\,k}_-$, we have
\begin{description}
\item[i.] $b_{\mu,E}(q)\leq{\underline{L}}_{\mu}^q(E)={\underline{C}}_{\mu}^q(E)$.
\item[ii.] ${\overline{L}}_{\mu,E}(q)={\overline{C}}_{\mu}^q(E)={\Lambda}_{\mu,E}(q)$.
\end{description}
\item For any $q\in\mathbb{R}^{*\,k}_+$, we have
$$
{\overline{L}}_{\mu,E}(q)\leq{\overline{C}}_{\mu}^q(E)\leq{\Lambda}_{\mu,E}(q).
$$
\end{enumerate}
\end{theorem}
{\bf Proof.} {\bf 1.} Using Besicovitch covering theorem we get
$$
\mathcal{T}_{\mu,\delta}^q(E)\leq\,C\mathcal{S}_{\mu,\delta}^q(E),
$$
with some constant $C$ fixed. So as {\bf 1.} is proved.\\
{\bf 2.} We firstly prove that
$$
{\underline{L}}_{\mu}^q(E)\geq{\underline{C}}_{\mu}^q(E)\quad\hbox{and}\quad{\overline{L}}_{\mu}^q(E)\geq{\overline{C}}_{\mu}^q(E).
$$
Indeed, let $\bigl(B(x_i,\delta)\bigr)_i$ be a centered $\delta$-packing of $E$ and $\bigl(B(y_i,\frac{\delta}{2})\bigr)$ be a centered $\frac{\delta}{2}$-covering of $E$. Consider for each $i$, the integer $k_i$ such that $x_i\in B(y_{k_i},\frac{\delta}{2})$. It is straightforward that for $i\not=j$ we have $k_i\not=k_j$. Consequently, for $q\in\mathbb{R}^{*\,k}_-$, there holds that
$$
\matrix{\displaystyle\sum_i\bigl(\mu(B(x_i,\delta))\bigr)^q
&=&\displaystyle\sum_i\Biggl(\displaystyle\frac{\mu(B(x_i,\delta))}{\mu(B(y_{k_i},\delta/2))}\Biggr)^q
\bigl(\mu(B(y_{k_i},\frac{\delta}{2}))\bigr)^q\hfill\cr\medskip
&\leq&\displaystyle\sum_i\bigl(\mu(B(y_i,\frac{\delta}{2}))\bigr)^q.\hfill}
$$
Which means that
$$
\mathcal{S}_{\mu,\delta}^q(E)\leq\mathcal{T}_{\mu,\frac{\delta}{2}}^q(E)
$$
and thus, for any $q\in\mathbb{R}^{*\,k}_-$,
$$
{\underline{L}}_{\mu}^q(E)\geq{\underline{C}}_{\mu}^q(E)\quad\hbox{and}\quad{\overline{L}}_{\mu}^q(E)\geq{\overline{C}}_{\mu}^q(E)
$$
Using the assertion {\bf 1.}, we obtain the equalities
$$
{\underline{L}}_{\mu}^q(E)={\underline{C}}_{\mu}^q(E)\quad\hbox{and}\quad{\overline{L}}_{\mu}^q(E)={\overline{C}}_{\mu}^q(E)
$$
for all $q\in\mathbb{R}^{*\,k}_-$. Therefore, to prove {\bf 2.i.}, it remains to prove the inequality of the left hand side. So, let $t>{\underline{L}}_\mu^q(E)$ and $F\subseteq E$. Consider next a sequence $(\delta_n)_n\subseteq]0,1[$ to be $\downarrow0$, and satisfying
$$
t>\displaystyle\frac{\log(\mathcal{T}_{\mu,\delta_n}^q(E))}{-\log\delta_n},\quad\forall\,n\in\mathbb{N}.
$$
This means that for each $n\in\mathbb{N}$, there exists a centered $\delta_n$-covering $\bigl(B(x_{ni},\delta_n)\bigr)_i$ of $E$ such that
$$
\displaystyle\sum_i\bigl(\mu(B(x_{ni},\delta_n))\bigr)^q<\delta_n^{-t}.
$$
There balls may be considered to be intersecting the set $F$. Next, for each $i$, choose an element $y_i\in\,B(x_{ni},\delta_n)\cap\,F$.
This results on a centered $2\delta_n$-covering $\bigl(B(y_i,2\delta_n)\bigr)_i$ of $F$. Therefore,
$$
\matrix{{\overline\mathcal{H}}_{\mu,2\delta_n}^{q,t}(F)
&\leq&\displaystyle\sum_i\bigl(\mu(B(x_{ni},\delta_n))\bigr)^q(4\delta_n)^t\hfill\cr\medskip
&=&4^t\displaystyle\sum_i\Biggl(\displaystyle\frac{\mu(B(y_i,2\delta_n))}{\mu(B(x_{ni},\delta_n))}\Biggr)^q
\bigl(\mu(B(x_{ni},\delta_n))\bigr)^q\delta_n^t\hfill\cr\medskip
&\leq&4^t\displaystyle\sum_i\bigl(\mu(B(x_{ni},\delta_n))\bigr)^q\delta_n^t\hfill\cr\medskip
&\leq&4^t\delta_n^{-t}\delta_n^t=4^t.\hfill}
$$
Hence,
$$
{\overline\mathcal{H}}_{\mu}^{q,t}(F)\leq 4^t,\quad\forall\,F\subseteq\,E,\;\;t>{\underline{L}}_\mu^q(E).
$$
So that,
$$
\mathcal{H}_{\mu}^{q,t}(E)\leq4^t<\infty,\quad\forall\,t>{\underline{L}}_\mu^q(E).
$$
Consequently,
$$
b_{\mu,E}(q)\leq t,\quad\forall\,t>{\underline{L}}_\mu^q(E)\,\Rightarrow\,b_{\mu,E}(q)\leq{\underline{L}}_\mu^q(E).
$$
We now prove the remaining par of {\bf 2.ii.} We will prove firstly that
\begin{equation}\label{bouligminkowski1}
{\overline{C}}_\mu^q(E)\leq\Lambda_{\mu,E}(q),\forall\,q\in\mathbb{R}^k.
\end{equation}
This is of course obvious when the right hand term is infinite. So, without loss of the generality, we assume that it is finite. Denote $t=\Lambda_{\mu,E}(q)$, and consider $\varepsilon>0$ and $0<\delta_\varepsilon<1$ be such that ${\overline\mathcal{P}}_{\mu,\delta}^{q,t+\varepsilon}(E)<1$ for all $0<\delta<\delta_\varepsilon$. This is possible because of the fact that ${\overline\mathcal{P}}_{\mu}^{q,t+\varepsilon}(E)=\displaystyle\lim_{\delta\downarrow0}{\overline\mathcal{P}}_{\mu,\delta}^{q,t+\varepsilon}(E)=0$. Consequently, for a centered $\delta$-packing $\bigl(B(x_i,\delta)\bigr)_i$ of $E$, we obatin
$$
\matrix{\displaystyle\sum_i\bigl(\mu(B(x_i,\delta))\bigr)^q
&=&(2\delta)^{-(t+\varepsilon)}\displaystyle\sum_i\bigl(\mu(B(x_i,\delta))\bigr)^q(2\delta)^{t+\varepsilon}\hfill\cr\medskip
&\leq&(2\delta)^{-(t+\varepsilon)}{\overline\mathcal{P}}_{\mu,\delta}^{q,t+\varepsilon}(E)\hfill\cr\medskip
&\leq&(2\delta)^{-(t+\varepsilon)}.\hfill}
$$
Hence,
$\mathcal{S}_{\mu,\delta}^q(E)\leq(2\delta)^{-(t+\varepsilon)}$ and consequently, equation (\ref{bouligandminkowskimixte}) holds. We now prove the converse
$$
\Lambda_{\mu,E}(q)\leq{\overline{C}}_\mu^q(E),\forall\,q\in\mathbb{R}^{*\,k}_-.
$$
Let $t=\Lambda_{\mu,E}(q)$, $\varepsilon>0$ and $0<\delta_0<1$. It holds that
$$
\infty={\overline\mathcal{P}}_{\mu}^{q,t-\varepsilon/2}(E)\leq{\overline\mathcal{P}}_{\mu,\delta_0}^{q,t-\varepsilon/2}(E).
$$
This means that there exists a centered $\delta_0$-packing $\bigl(B(x_i,r_i)\bigr)_i$ of $E$ such that
$$
1<\displaystyle\sum_i\bigl(\mu(B(x_i,r_i))\bigr)^q(2r_i)^{t-\varepsilon/2}.
$$
Next, denote for $n\in\mathbb{N}$,
$$
I_n=\{\ i\in\mathbb{N}\,;\;\;\displaystyle\frac{\delta_0}{2^{n+1}}\leq\,r_i<\displaystyle\frac{\delta_0}{2^n}\,\}\quad\hbox{and}\quad
\nu_n=\displaystyle\sum_{i\in\,I_n}\bigl(\mu(B(x_i,r_i))\bigr)^q.
$$
A straightforward computation yields that
$$
C\displaystyle\sup_m\Biggl(\nu_m(\frac{\delta_0}{2^m})^{t-\varepsilon}\Biggr)>1
$$
for an appropriate constant $C>0$ depending only on $t$ and $\varepsilon$. Consequently, for $N\in\mathbb{N}$ such that
$1<C\nu_N(\frac{\delta_0}{2^N})^{t-\varepsilon}$ and $\delta=\frac{\delta_0}{2^{N+1}}$, the set $\bigl(B(x_i,\delta)\bigr)_i$ forms a centered $\delta$-packing de $E$. Observing that $q\in\mathbb{R}^{*\,k}_-$, it results that
$$
\mathcal{S}_{\mu,\delta}^q(E)\geq\,C^{-1}\delta^{-(t-\varepsilon)}.
$$
Consequently, ${\overline{C}}_\mu^q(E)\geq\Lambda_{\mu,E}(q)$ for all $q\in\mathbb{R}^{*\,k}_-$.\\
{\bf 3.} It follows from {\bf 1.} and equation (\ref{bouligminkowski1}).

Next we need to introduce the following quantities which will be useful later. Let $\mu=(\mu_1,\,\mu_2,\dots,\mu_k)$ be a vector valued measure composed of probability measures on $\mathbb{R}^d$. For $j=1,\,2,\,\dots,\,k$, $a>1$ and $E\subseteq\,Support(\mu)$, denote
$$
T_a^j(E)=\displaystyle\limsup_{r\downarrow0}\Bigl(\displaystyle\sup_{x\in\,E}\frac{\mu_j\bigl(B(x,ar)\bigr)}{\mu_j\bigl(B(x,r)\bigr)}\Bigr)
$$
and for $x\in\,Support(\mu)$, $T_a^j(x)=T_a^j(\{x\})$. Denote also
$$
P_0(\mathbb{R}^d,E)=\{\,\mu\,;\;\;\exists\,a\,>\,1\,;\;\;\forall\,x\in\,E,\;\;T_a^j(x)<\infty,\;\;\forall\,j\,\},
$$
$$
P_1(\mathbb{R}^d,E)=\{\,\mu\,;\;\;\exists\,a\,>\,1\,;\;T_a^j(E)<\infty,\;\;\forall\,j\,\},
$$
$$
P_0(\mathbb{R}^d)=P_0(\mathbb{R}^d,Support(\mu))\qquad\hbox{and}\qquad\,P_1(\mathbb{R}^d)=P_1(\mathbb{R}^d,Support(\mu)).
$$
\begin{theorem}\label{anouartheorem2}
\begin{enumerate}
\item For $\mu\in\,P_0(\mathbb{R}^d)$ and $q\in\mathbb{R}^{*\,k}_+$, there holds that
$$
b_{\mu,E}(q)\leq{\overline{L}}_{\mu}^q(E).
$$
\item For $\mu\in\,P_1(\mathbb{R}^d)$ and $q\in\mathbb{R}^{*\,k}_+$, there holds that
\begin{description}
\item[i.] ${\underline{L}}_{\mu}^q(E)={\underline{C}}_{\mu}^q(E)$.
\item[ii.] ${\overline{L}}_{\mu,E}(q)={\overline{C}}_{\mu}^q(E)={\Lambda}_{\mu,E}(q)$.
\end{description}
\end{enumerate}
\end{theorem}
{\bf Proof.} {\bf 1.} The vector valued measure $\mu\in\,P_0(\mathbb{R}^d)$ yields that
$$
E=\displaystyle\bigcup_{m\in\mathbb{N}}E_m
$$
where
$$
E_m=\{\,x\in\,E\,;\;\;\displaystyle\frac{\mu_j(B(x_i,4r))}{\mu_j(B(x_i,r))}<m\,,\;0<r<\frac{1}{m},\,\;\;\forall\,j\,\}.
$$
Next, remark that for $t>{\overline{L}}_\mu^q(E)$ and $F\subseteq\,E_m$, there exists a sequence $(\delta_n)_n\in]0,1[\downarrow0$ for which
$$
t<\displaystyle\frac{\log(\mathcal{T}_{\mu,\delta_n}^q(F))}{-\log\delta_n},\qquad\forall\,n\in\mathbb{N}.
$$
Therefore, there exists a centered $\delta_n$-covering $(B(x_{ni},\delta_n))_i$ of $F$ satisfying
$$
\displaystyle\sum_i\bigl(\mu(B(x_{ni},\delta_n))\bigr)^q<\delta_n^{-t}.
$$
Let next $y_{ni}\in\,B(x_{ni},\delta_n)$. Then, $(B(x_{ni},2\delta_n))_i$ is a centered $2\delta_n$-covering of $F$. Hence,
$$
\matrix{{\overline\mathcal{H}}_{\mu,2\delta_n}^{q,t}(F)&\leq&\displaystyle\sum_i\bigl(\mu(B(y_{ni},2\delta_n))\bigr)^q(4\delta_n)^t\hfill\cr\medskip
&\leq&4^t\displaystyle\sum_i\Biggl(\displaystyle\frac{\mu(B(y_{ni},2\delta_n))}{\mu(B(x_{ni},\delta_n))}\Biggr)^q
\bigl(\mu(B(x_{ni},\delta_n))\bigr)^q\delta_n^t\hfill\cr\medskip
&\leq&4^t\displaystyle\sum_i\Biggl(\displaystyle\frac{\mu(B(x_{ni},4\delta_n))}{\mu(B(x_{ni},\delta_n))}\Biggr)^q
\bigl(\mu(B(x_{ni},\delta_n))\bigr)^q\delta_n^t\hfill\cr\medskip
&\leq&4^tm^{|q|}\displaystyle\sum_i\bigl(\mu(B(x_{ni},\delta_n))\bigr)^q\delta_n^t\hfill\cr\medskip
&\leq&4^tm^{|q|}\hfill}
$$
where $|q|=q_1+q_2+\dots+q_k$. Thus,
$$
{\overline\mathcal{H}}_{\mu}^{q,t}(F)\leq4^tm^{|q|},\quad\forall\,m,\;\;\hbox{and}\;\,F\subseteq\,E_m.
$$
Which means that
$$
\mathcal{H}_{\mu}^{q,t}(E_m)\leq4^tm^{|q|}<\infty,\quad\forall\,m,\;\;\hbox{and}\;\,t>{\underline{L}}_\mu^q(E).
$$
Consequently,
$$
b_{\mu,E_m}(q)\leq\,t,\quad\forall\,m,\;\;\hbox{and}\;\,t>{\underline{L}}_\mu^q(E).
$$
Using the $\sigma$-stability of $b_{\mu,.}(q)$ (See Proposition \ref{anouarproposition2}. c.), it results that
$$
b_{\mu,E}(q)\leq\,t,\quad\forall\,t>{\underline{L}}_\mu^q(E)\;\Rightarrow\;b_{\mu,E}(q)\leq{\underline{L}}_\mu^q(E).
$$
{\bf 2. i.} From Theorem \ref{anouartheorem1}, it remains to prove that ${\underline{L}}_{\mu}^q(E)\geq{\underline{C}}_{\mu}^q(E)$. Let $C,\delta_0>0$ such that
$$
\displaystyle\sup_{x\in\,E}\displaystyle\frac{\mu_j(B(x,4r))}{\mu(B(x,r))}<C,\quad\forall\,0<r<\delta_0,\quad\forall\,j=1,2,\dots,k.
$$
Let next $0<\delta<\delta_0$, $(B(x_i,\delta))_i$ be a centered packing and $(B(y_i,\delta/2))_i$ be a centered covering of $E$. For each $i\in\mathbb{N}$, denote $k_i$ the unique integer such that $x_i\in\,B(y_{k_i},\delta/2)$. It holds that
$$
\matrix{\mathcal{S}_{\mu,\delta}^q(E)&\leq&\displaystyle\sum_i\bigl(\mu(B(x_{i},\delta))\bigr)^q\hfill\cr\medskip
&=&\displaystyle\sum_i\Biggl(\displaystyle\frac{\mu(B(x_{i},\delta))}{\mu(B(y_{k_i},\delta/2))}\Biggr)^q\bigl(\mu(B(y_{k_i},\delta/2))\bigr)^q\hfill\cr\medskip
&\leq&\displaystyle\sum_i\Biggl(\displaystyle\frac{\mu(B(y_{k_i},2\delta))}{\mu(B(y_{k_i},\delta/2))}\Biggr)^q\bigl(\mu(B(y_{k_i},\delta/2))\bigr)^q\hfill\cr\medskip
&\leq&C^{|q|}\displaystyle\sum_i\bigl(\mu(B(y_{k_i},\delta/2))\bigr)^q\hfill\cr\medskip
&\leq&C^{|q|}\displaystyle\sum_i\bigl(\mu(B(y_{i},\delta/2))\bigr)^q.\hfill}
$$
This yields that
$$
\mathcal{S}_{\mu,\delta}^q(E)\leq\,C^{|q|}\mathcal{T}_{\mu,\delta/2}^q(E).
$$
Consequently,
$$
{\underline{C}}_{\mu}^q(E)\leq{\underline{L}}_{\mu}^q(E)\quad\hbox{and}\quad{\overline{C}}_{\mu}^q(E)\leq{\overline{L}}_{\mu}^q(E)
$$
Using Theorem \ref{anouartheorem1}, {1.}, we obtain the equalities.
\begin{equation}\label{bouligminkowski2}
{\underline{C}}_{\mu}^q(E)={\underline{L}}_{\mu}^q(E)\quad\hbox{and}\quad{\overline{C}}_{\mu}^q(E)={\overline{L}}_{\mu}^q(E),
\quad\forall\,q\in\mathbb{R}^{*\,k}_+,\;\;\mu\in\,P_1(\mathbb{R}^d)
\end{equation}
{\bf ii.} Using equations (\ref{bouligminkowski1}) and (\ref{bouligminkowski2}), it remains to prove that
$$
{\overline{C}}_\mu^q(E)\geq\Lambda_{\mu,E}(q).
$$
The vector measure $\mu$ lies in $P_1{\mathbb{R}^d,E}$. So that, there exists as above $C>0$, and $0<r_0<1$ such that
$$
\displaystyle\frac{\mu_j(B(x,2r))}{\mu_j(B(x,r))}\leq\,C,\;\;\forall\,0<r<\delta_0,\;\;\,x\in\,E,\;\;\hbox{and}\;\;j=1,2,\dots,k.
$$
Denote $t=\Lambda_{\mu,E}(q)$, $\varepsilon>0$  and $0<\delta_0<r_0$. Then ${\overline\mathcal{P}}_\mu^{q,t-\varepsilon/2}(E)=\infty$. Which means that there exists a centered $\delta_0$-packing of the set $E$, $(B(x_i,r_i))_i$ for which
$$
1<\displaystyle\sum_i\Bigl(\mu(B(x_i,r_i))\Bigr)^q(2r_i)^{t-\varepsilon/2}.
$$
By considering the set $I_N$, $\nu_N$ and $\delta$ as above, we obtain
$$
\matrix{{\cal S}_{\mu,\delta}^q(E)&\geq&\displaystyle\sum_{i\in\,I_N}\Bigl(\mu(B(x_i,\delta))\Bigr)^q\hfill\cr\medskip
&\geq&\beta^{-q}\displaystyle\sum_{i\in\,I_N}\Biggl(\frac{\mu(B(x_i,\delta_0/2^{N+1}))}
{\mu(B(x_i,\delta_0/2^{N}))}\Biggr)^q\Bigl(\mu(B(x_i,r_i))\Bigr)^q\hfill\cr\medskip
&\geq&C^{-|q|}\displaystyle\sum_{i\in\,I_N}\Bigl(\mu(B(x_i,r_i))\Bigr)^q\hfill\cr\medskip
&\geq&C^{-|q|}\nu_N>\beta^{-q}C^{-1}(\displaystyle\frac{\delta_0}{2^N})^{-(t-\varepsilon)}.\hfill}
$$
Hence,
$$
{\overline{C}}_\mu^q(E)\geq\Lambda_{\mu,E}(q).
$$
We now recall re-introduce the mixed multifrcatal generalization of the $L^q$-dimensions called also Renyi dimensions based on integral representations. See \cite{olsen3} for more details and other results. For $q\in\mathbb{R}^{*,k}$, $\mu=(\mu_1,\,\mu_2,\,\dots,\,\mu_k)$ and $\delta>0$, we set
$$
I_{\mu,\delta}^q=\displaystyle\int_{S_{\mu}}\Bigl(\mu(B(t,\delta))\Bigr)^qd\mu(t),
$$
where, in this case,
$$
S_{\mu}=Support(\mu_1)\times\,Support(\mu_2)\times\,\dots\,\times\,Support(\mu_k),
$$
$$
\Bigl(\mu(B(t,\delta))\Bigr)^qd\mu(t)=\Bigl(\mu_1(B(t_1,\delta))\Bigr)^{q_1}\,\Bigl(\mu_2(B(t_2,\delta))\Bigr)^{q_2}\dots\Bigl(\mu_k(B(t_k,\delta))\Bigr)^{q_k}
$$
and
$$
d\mu(t)=d\mu_1(t_1)\,d\mu_2(t_2)\,\dots\,d\mu_k(t_k).
$$
The mixed multifractal generalizations of the Renyi dimensions are
$$
{\overline{I}}_\mu^q=\displaystyle\limsup_{\delta\downarrow0}\displaystyle\frac{\log\,I_{\mu,\delta}^q}{-\log\delta},\quad\hbox{and}\quad
{\underline{I}}_\mu^q=\displaystyle\liminf_{\delta\downarrow0}\displaystyle\frac{\log\,I_{\mu,\delta}^q}{-\log\delta}.
$$
We now propose to relate these dimensions to the quantities ${\underline{C}}_\mu^{q}$, ${\overline{C}}_\mu^{q}$, ${\underline{L}}_\mu^{q}$, ${\overline{L}}_\mu^{q}$ introduce previously.
\begin{proposition}\label{anouarproposition5} The following results hold.
\begin{description}
\item[a.] $\forall\,q\in\mathbb{R}^{*,\,k}_-$,
$$
\underline{C}_{\mu}^{q+\1}(Support(\mu))\geq\underline{I}_{\mu}^q\quad\,and\,\quad\overline{C}_{\mu}^{q+\1}(Support(\mu))\geq\overline{I}_{\mu}^q.
$$
\item[b.] $\forall\,q\in\mathbb{R}^{*,\,k}_+$,
$$
{\underline{C}}_\mu^{q+\1}(supp(\mu))\leq{\underline{I}_\mu^q}\qquad\,and\qquad\,{\overline{C}}_\mu^{q+\1}(supp(\mu))\leq{\overline{I}_\mu^q}.
$$
\item[c.] $\forall\,q\in\mathbb{R}^{*,\,k}$, $\mu\in\,P_1(\mathbb{R}^d)$,
$$
{\underline{C}}_\mu^{q+\1}(supp(\mu))=\underline{I}_\mu^q\qquad\,and\qquad\,{\overline{C}}_\mu^{q+\1}(supp(\mu))=\overline{I}_\mu^q.
$$
\item[d.] $\forall\,q\in\mathbb{R}^{*,\,k}_-$,
$$
\underline{I}_\mu^q\leq\underline{L}_\mu^{q+\1}(supp(\mu))\qquad\,and\qquad\,\overline{I}_\mu^q\leq\overline{L}_\mu^{q+\1}(supp(\mu)).
$$
\end{description}
\end{proposition}
{\bf Proof.} {\bf a.} For $\delta>0$, let $\Bigl(B(x_i,\delta)\Bigr)_i$ be a centered $\delta$-covering of $Support(\mu)$ and let next $\Bigl(B(x_{ij},\delta)\Bigr)_j$, ${1\leq\,i\leq\xi}$ the $\xi$ sets defined in Besicovitch covering theorem. It holds that
$$
\matrix{\displaystyle\sum_{i,j}\Bigl(\mu(B(x_{ij},\delta))\Bigr)^{q+\1}
&=&\displaystyle\sum_{i,j}\Bigl(\mu(B(x_{ij},\delta))\Bigr)^q\displaystyle\int_{B(x_{ij},\delta)^k}d\mu(t)\hfill\cr\medskip
&\geq&\displaystyle\sum_{i,j}\displaystyle\int_{B(x_{ij},\delta)^k}\Bigl(\mu(B(t,2\delta))\Bigr)^qd\mu(t)\hfill\cr\medskip
&\geq&\displaystyle\int_{S_{\mu}}\Bigl(\mu(B(t,2\delta))\Bigr)^qd\mu(t).\hfill}
$$
As a results,
$$
\xi\mathcal{S}_{\mu,\delta}^{q+\1}(Support(\mu))\geq\,I_{\mu,2\delta}^q.
$$
Which implies that
$$
\underline{C}_{\mu}^{q+\1}(Support(\mu))\geq\underline{I}_{\mu}^q\quad\hbox{and}\quad\overline{C}_{\mu}^{q+\1}(Support(\mu))\geq\overline{I}_{\mu}^q.
$$
{\bf b.} Let $\delta>0$  and $\Bigl(B(x_i,\delta)\Bigr)_i$ a centered $\delta$-packing of $Support(\mu)$. It holds that
$$
\matrix{\displaystyle\sum_i\Bigl(\mu(B(x_i,\delta))\Bigr)^{q+\1}
&=&\displaystyle\sum_i\Bigl(\mu(B(x_i,\delta))\Bigr)^q\displaystyle\int_{B(x_i,\delta)^k}d\mu(t)\hfill\cr\medskip
&\leq&\displaystyle\sum_i\displaystyle\int_{B(x_i,\delta)^k}\Bigl(\mu(B(t,2\delta))\Bigr)^qd\mu(t)\hfill\cr\medskip
&\leq&\displaystyle\int_{S_{\mu}}\Bigl(\mu(B(t,2\delta))\Bigr)^qd\mu(t).\hfill}
$$
Therefore,
$$
\mathcal{S}_{\mu,\delta}^{q+\1}(Support(\mu))\leq\,I_{\mu,2\delta}^q
$$
and thus,
$$
\underline{C}_{\mu}^{q+\1}(Support(\mu))\leq\underline{I}_{\mu}^q\quad\hbox{and}\quad\overline{C}_{\mu}^{q+\1}(Support(\mu))\leq\overline{I}_{\mu}^q.
$$
{\bf c.} Assume firstly that $q\in\mathbb{R}^{*,\,k}_-$. Observing assertion {\bf a.}, it suffices to prove that
$$
\underline{C}_{\mu}^{q+\1}(Support(\mu))\leq\underline{I}_{\mu}^q\quad\,and\,\quad\overline{C}_{\mu}^{q+\1}(Support(\mu))\leq\overline{I}_{\mu}^q.
$$
Since the measure $\mu\in\,P_1(\mathbb{R}^d)$, there exists a constant $C>0$ and $r_0>0$ such that
$$
\displaystyle\frac{\mu_j(B(x,2r))}{\mu_j(B(x,r))}<C\,;\;\;\forall\,x\in\,Support(\mu),\;\;0<r<r_0,\;\;j=1,2,\dots\,k.
$$
Next, consider for $0<\delta<r_0$ a centered $\delta$-packing $(B(x_i,\delta))_i$ of $Support(\mu)$. It holds that
$$
\matrix{\displaystyle\sum_i\Bigl(\mu(B(x_i,\delta))\Bigr)^{q+\1}
&=&\displaystyle\sum_i\Bigl(\mu(B(x_i,\delta))\Bigr)^q\displaystyle\int_{B(x_i,\delta)^k}d\mu(t)\hfill\cr\medskip
&\leq&C^{-2|q|}\displaystyle\sum_i\displaystyle\int_{B(x_i,\delta)^k}\Bigl(\mu(B(t,2\delta))\Bigr)^qd\mu(t)\hfill\cr\medskip
&\leq&C^{-2|q|}\displaystyle\int_{S_{\mu}}\Bigl(\mu(B(t,2\delta))\Bigr)^qd\mu(t).\hfill}
$$
Consequently,
$$
\mathcal{S}_{\mu,\delta}^{q+\1}(Support(\mu))\leq\,C^{-2|q|}I_{\mu,2\delta}^q.
$$
Hence,
$$
\underline{C}_{\mu}^{q+\1}(Support(\mu))\leq\underline{I}_{\mu}^q\quad\hbox{and}\quad\overline{C}_{\mu}^{q+\1}(Support(\mu))\leq\overline{I}_{\mu}^q.
$$
So the equality for $q\in\mathbb{R}^{*,\,k}_-$. \\
Assume now that $q\in\mathbb{R}^{*,\,k}_+$. Observing assertion {\bf b.}, it remains to prove that
$$
{\underline{C}}_\mu^{q+\1}(supp(\mu))\geq{\underline{I}_\mu^q}\qquad\,and\qquad\,{\overline{C}}_\mu^{q+\1}(supp(\mu))\geq{\overline{I}_\mu^q}.
$$
To do so, we use the fact that $\mu\in\,P_1(\mathbb{R}^d)$, which means that there exists $C>0$ and $r_0>0$ satisfying
$$
\displaystyle\frac{\mu_j(B(x,2r))}{\mu_j(B(x,r))}<C,\;\;\forall\,x\in\,Support(\mu),\;\;0<r<r_0,\;\;j=1,2,\dots,k.
$$
Let next, $0<\delta<r_0$, $(B(x_i,\delta))_i$ a centered $\delta$-covering of $Support(\mu)$ and as previously, $\Bigl(B(x_{ij},\delta)\Bigr)_j$, ${1\leq\,i\leq\xi}$ the $\xi$ sets defined in Besicovitch covering theorem. We have
$$
\matrix{\displaystyle\sum_{i,j}\Bigl(\mu(B(x_{ij},\delta))\Bigr)^{q+\1}
&=&\displaystyle\sum_{i,j}\Bigl(\mu(B(x_{ij},\delta))\Bigr)^q\displaystyle\int_{B(x_{ij},\delta)^k}d\mu(t)\hfill\cr\medskip
&\geq&C^{-2|q|}\displaystyle\sum_{i,j}\displaystyle\int_{B(x_{ij},\delta)^k}\Bigl(\mu(B(t,2\delta))\Bigr)^qd\mu(t)\hfill\cr\medskip
&\geq&C^{-2|q|}\displaystyle\int_{S_{\mu}}\Bigl(\mu(B(t,2\delta))\Bigr)^qd\mu(t).\hfill}
$$
Hence,
$$
\xi\mathcal{S}_{\mu,\delta}^{q+\1}(Support(\mu))\geq\,C^{-2|q|}I_{\mu,2\delta}^q.
$$
Consequently,
$$
\underline{C}_{\mu}^{q+\1}(Support(\mu))\geq\underline{I}_{\mu}^q\quad\hbox{and}\quad\overline{C}_{\mu}^{q+\1}(Support(\mu))\geq\overline{I}_{\mu}^q.
$$
Hence, the equality for $q\in\mathbb{R}^{*,\,k}_+$.\\
{\bf d.} Let $\delta>0$ and $\Bigl(B(x_i,\delta)\Bigr)_i$ be a centered $\delta$-covering of $Support(\mu)$. We have
$$
\matrix{\displaystyle\sum_i\Bigl(\mu(B(x_i,\delta))\Bigr)^{q+\1}
&=&\displaystyle\sum_i\Bigl(\mu(B(x_i,\delta))\Bigr)^q\displaystyle\int_{B(x_i,\delta)^k}d\mu(t)\hfill\cr\medskip
&\geq&\displaystyle\sum_i\displaystyle\int_{B(x_i,\delta)^k}\Bigl(\mu(B(t,2\delta))\Bigr)^qd\mu(t)\hfill\cr\medskip
&\geq&\displaystyle\int_{S_{\mu}}\Bigl(\mu(B(t,2\delta))\Bigr)^qd\mu(t).\hfill}
$$
As a result,
$$
\mathcal{T}_{\mu,\delta}^{q+1}(Spport(\mu))\geq\,I_{\mu,2\delta}^q.
$$
Consequently,
$$
\underline{L}_{\mu}^{q+\1}(supp(\mu))\geq\underline{I}_{\mu}^q\quad\hbox{and}\quad\overline{L}_{\mu}^{q+\1}(supp(\mu))\geq\overline{I}_{\mu}^q.
$$
\section{A mixed multifractal formalism for vector valued measures}
Let $\mu=(\mu_1,\,\mu_2,\,\dots,\,\mu_k)$ be a vector valued probability measure on $\mathbb{R}^d$. For $x\in\mathbb{R}^d$ and $j=1,2,\dots,k$, we denote
$$
{\underline\alpha}_{\mu_j}(x)=\displaystyle\liminf_{r\downarrow0}\displaystyle\frac{\log(\mu_j(B(x,r)))}{\log\,r}\;\;\hbox{and}\;\;
{\overline\alpha}_{\mu_j}(x)=\displaystyle\limsup_{r\downarrow0}\displaystyle\frac{\log(\mu_j(B(x,r)))}{\log\,r}
$$
respectively the local lower dimension and the local upper dimension of $\mu_j$ at the point $x$ and as usually the local dimension $\alpha_{\mu_j}(x)$ of $\mu_j$ at $x$ will be the common value when these are equal. Next for $\alpha=(\alpha_1,\,\alpha_2,\,\dots,\,\alpha_k)\in\mathbb{R}_+^k$, let
$$
\underline{X}_\alpha=\{\,x\in\,Support(\mu)\,;\,\,{\underline\alpha}_{\mu_j}(x)\geq\alpha_j\,,\forall\,j=1,2,\dots,k\,\},
$$
$$
\overline{X}^\alpha=\{\,x\in\,Support(\mu)\,;\,\,{\overline\alpha}_{\mu_j}(x)\leq\alpha_j\,,\forall\,j=1,2,\dots,k\,\}
$$
and
$$
X(\alpha)=\underline{X}_\alpha\cap\overline{X}^\alpha.
$$
The mixed multifractal spectrum of the vector valued measure $\mu$ is defined by
$$
\alpha\,\longmapsto\,dim\,X(\alpha)
$$
where $dim$ stands for the Hausdorff dimension.

In this section, we propose to compute such a spectrum for some cases of measures that resemble to the situation raised by Olsen in \cite{olsen1} but in the mixed case. This will permit to describe better the simultaneous behavior of finitely many measures. We intend precisely to compute the mixed spectrum based on the mixed multifractal generalizations of the Haudorff and packing dimensions $b_\mu$, $B_\mu$ and $\Lambda_\mu$. We start with the following technic results.
\begin{lemma}\label{spectrelemme1}
\begin{description}
\item[1.] $\forall\,\delta>0,\,t\in\mathbb{R}$ and $q\in\mathbb{R}^k_+$, $\alpha\in\mathbb{R}^k$ such that $\langle\alpha,q\rangle+t\geq0$, we have
\begin{description}
\item[i.] $\mathcal{H}^{\langle\alpha,q\rangle+t+k\delta}(\overline{X}^\alpha)\leq2^{\langle\alpha,q\rangle+k\delta}\mathcal{H}_\mu^{q,t}(\overline{X}^\alpha).$
\item[ii.]
$\mathcal{P}^{\langle\alpha,q\rangle+t+k\delta}({\overline{X}}^\alpha)\leq2^{\langle\alpha,q\rangle+k\delta}\mathcal{P}_\mu^{q,t}({\overline{X}}^\alpha).$
\end{description}
\item[2.] $\forall\,\delta>0,\,t\in\mathbb{R}$ and $q\in\mathbb{R}^k_-$, $\alpha\in\mathbb{R}^k$ such that $\langle\alpha,q\rangle+t\geq0$, we have
\begin{description}
\item[i.] $\mathcal{H}^{\langle\alpha,q\rangle+t+k\delta}({\underline{X}}_\alpha)\leq2^{\langle\alpha,q\rangle+k\delta}\mathcal{H}_\mu^{q,t}({\underline{X}}_\alpha).$
\item[ii.] $\mathcal{P}^{\langle\alpha,q\rangle+t+k\delta}({\underline{X}}_\alpha)\leq2^{\langle\alpha,q\rangle+k\delta}\mathcal{P}_\mu^{q,t}({\underline{X}}_\alpha).$
\end{description}
\end{description}
\end{lemma}
{\it Proof.} {\bf 1.} {$\mathbf{i.}$} We prove the first part. For $m\in\mathbb{N}^*$, consider the set
$$
\overline{X}_m^\alpha=\{\,x\in\overline{X}^\alpha;\,\,\displaystyle\frac{\log(\mu_j(B(x,r)))}{\log\,r}\leq\alpha_j+\displaystyle\frac{\delta}{q_j};\,\,
0<r<\displaystyle\frac{1}{m},\;\;1\leq\,j\leq\,k\,\}.
$$
Let next $0<\eta<\displaystyle\frac{1}{m}$ and $(B(x_i,r_i))_i$ a centered $\eta$-covering of ${\overline X}_m^\alpha$. It holds that
$$
(\mu(B(x_i,r_i)))^{q}\geq\,r_i^{\langle\alpha,q\rangle+k\delta}.
$$
Consequently,
$$
\mathcal{H}_\eta^{\langle\alpha,q\rangle+t+k\delta}(\overline{X}_m^\alpha)\leq\displaystyle\sum_i(2r_i)^{\langle\alpha,q\rangle+t+k\delta}
\leq2^{\langle\alpha,q\rangle+k\delta}\displaystyle\sum_i(\mu(B(x_i,r_i)))^q(2r_i)^t.
$$
Hence, $\forall\eta>0$, there holds that
$$
\mathcal{H}_\eta^{\langle\alpha,q\rangle+t+k\delta}(\overline{X}_m^\alpha)
\leq2^{\langle\alpha,q\rangle+k\delta}{\overline\mathcal{H}}_{\mu,\eta}^{q,t}(\overline{X}_m^\alpha).
$$
Which means that
$$
\mathcal{H}^{\langle\alpha,q\rangle+t+k\delta}(\overline{X}_m^\alpha)
\leq2^{\langle\alpha,q\rangle+k\delta}{\overline\mathcal{H}}_{\mu}^{q,t}(\overline{X}_m^\alpha)
\leq2^{\langle\alpha,q\rangle+k\delta}\mathcal{H}_{\mu}^{q,t}(\overline{X}_m^\alpha).
$$
Next, observing that $\overline{X}^\alpha=\displaystyle\bigcup_{m}\overline{X}_m^\alpha$, we obtain
$$
\mathcal{H}^{\langle\alpha,q\rangle+t+k\delta}(\overline{X}^\alpha)
\leq2^{\langle\alpha,q\rangle+k\delta}\mathcal{H}_{\mu}^{q,t}(\overline{X}^\alpha).
$$
{$\mathbf{ii.}$} For $q\in\mathbb{R}^{*,k}_+$ and $m\in\mathbb{N}^*$, consider the set $\overline{X}_m^\alpha$ defined previously and let $E\subseteq\overline{X}_m^\alpha$, $0<\eta<\displaystyle\frac{1}{m}$ and $\bigl(B(x_i,r_i)\bigr)_i$ a centered $\eta$-packing of $E$. We have
$$
\displaystyle\sum_i(2r_i)^{\langle\alpha,q\rangle+t+k\delta}\leq2^{\langle\alpha,q\rangle+k\delta}\displaystyle\sum_i(\mu(B(x_i,r_i)))^q(2r_i)^t
\leq2^{\langle\alpha,q\rangle+k\delta}{\overline\mathcal{P}}_{\mu,\eta}^{q,t}(E).
$$
Consequently, $\forall\,\eta>0$,
$$
{\overline\mathcal{P}}_\eta^{\langle\alpha,q\rangle+t+k\delta}(E)\leq2^{\langle\alpha,q\rangle+k\delta}{\overline\mathcal{P}}_{\mu,\eta}^{q,t}(E).
$$
Hence, $\forall\,E\subseteq\overline{X}_m^\alpha$,
$$
{\overline\mathcal{P}}^{\langle\alpha,q\rangle+t+k\delta}(E)\leq2^{\langle\alpha,q\rangle+k\delta}{\overline\mathcal{P}}_{\mu}^{q,t}(E).
$$
Let next, $(E_i)_i$ be a covering of $\overline{X}_m^\alpha$. Thus,
$$
\matrix{\mathcal{P}^{\langle\alpha,q\rangle+t+k\delta}(\overline{X}_m^\alpha)
&=&\mathcal{P}^{\langle\alpha,q\rangle+t+k\delta}\Biggl(\displaystyle\bigcup_i(\overline{X}_m^\alpha\cap\,E_i)\Biggr)\hfill\cr\medskip
&=&\displaystyle\sum_i\mathcal{P}^{\langle\alpha,q\rangle+t+k\delta}\Bigl(\overline{X}_m^\alpha\cap\,E_i\Bigr)\hfill\cr\medskip
&\leq&\displaystyle\sum_i{\overline\mathcal{P}}^{\langle\alpha,q\rangle+t+k\delta}\Bigl(\overline{X}_m^\alpha\cap\,E_i\Bigr)\hfill\cr\medskip
&\leq&2^{\langle\alpha,q\rangle+k\delta}\displaystyle\sum_i{\overline\mathcal{P}}_\mu^{q,t}\Bigl(\overline{X}_m^\alpha\cap\,E_i\Bigr)\hfill\cr\medskip
&\leq&2^{\langle\alpha,q\rangle+k\delta}\displaystyle\sum_i{\overline\mathcal{P}}_\mu^{q,t}(E_i).\hfill}
$$
Hence, $\forall,m$,
$$
\mathcal{P}^{\langle\alpha,q\rangle+t+k\delta}(\overline{X}_m^\alpha)\leq2^{\langle\alpha,q\rangle+k\delta}\mathcal{P}_\mu^{q,t}(\overline{X}_m^\alpha).
$$
Consequently,
$$
\mathcal{P}^{\langle\alpha,q\rangle+t+k\delta}(\overline{X}^\alpha)\leq2^{\langle\alpha,q\rangle+k\delta}\mathcal{P}_\mu^{q,t}(\overline{X}^\alpha).
$$
{\bf 2.} {$\mathbf{i.}$} and {$\mathbf{ii.}$} follow similar arguments and techniques as previously.
\begin{proposition}\label{spectreproposition1}
Let $\alpha\in\mathbb{R}^k_+$ and $q\in\mathbb{R}^k$. The following assertions hold.
\begin{description}
\item[a.] Whenever $\langle\alpha,q\rangle+b_\mu(q)\geq0$, we have
\begin{description}
\item[i.] $dim{\overline{X}}^\alpha\leq\langle\alpha,q\rangle+b_\mu(q)$,\quad$\forall\,q\mathbb{R}^k_+$.
\item[ii.] $dim{\underline{X}}_\alpha\leq\langle\alpha,q\rangle+b_\mu(q)$\quad$\forall\,q\mathbb{R}^k_-$.
\end{description}
\item[b.] Whenever $\langle\alpha,q\rangle+B_\mu(q)\geq0$, we have
\begin{description}
\item[i.] $Dim{\overline{X}}^\alpha\leq\langle\alpha,q\rangle+B_\mu(q)$,\quad$\forall\,q\mathbb{R}^k_+$.
\item[ii.] $Dim{\underline{X}}_\alpha\leq\langle\alpha,q\rangle+B_\mu(q)$\quad$\forall\,q\mathbb{R}^k_-$.
\end{description}
\end{description}
\end{proposition}
{\bf Proof.} {\bf a. i.} It follows from Lemma \ref{spectrelemme1}, assertion {\bf 1. i.},
$$
\mathcal{H}^{\langle\alpha,q\rangle+t+k\delta}({\overline{X}}^\alpha)=0,\qquad\forall\,t>b_\mu(q),\,\,\delta>0.
$$
Consequently,
$$
dim{\overline{X}}^\alpha\leq\langle\alpha,q\rangle+t+k\delta,\qquad\forall\,t>b_\mu(q),\;\;\delta>0.
$$
Hence,
$$
dim\,\overline{X}^\alpha\leq\langle\alpha,q\rangle+b_\mu(q).
$$
{\bf a. ii.} It follows from Lemma \ref{spectrelemme1}, assertion {\bf 2. i.}, as previously, that
$$
\mathcal{H}^{\langle\alpha,q\rangle+t+k\delta}(\underline{X}^\alpha)=0,\qquad\forall\,t>b_\mu(q),\,\,\delta>0.
$$
Hence,
$$
dim{\underline{X}}_\alpha\leq\langle\alpha,q\rangle+t+k\delta,\qquad\forall\,t>b_\mu(q),\;\;\delta>0
$$
and finally,
$$
dim{\underline{X}}_\alpha\leq\langle\alpha,q\rangle+b_\mu(q).
$$
{\bf b. i.} observing Lemma \ref{spectrelemme1}, assertion {\bf 1. ii.}, we obtain
$$
\mathcal{P}^{\langle\alpha,q\rangle+t+k\delta}(\overline{X}^\alpha),\qquad\forall\,t>B_\mu(q),\;\;\delta>0.
$$
Consequently,
$$
Dim\,\overline{X}^\alpha\leq\langle\alpha,q\rangle+t+k\delta,\qquad\forall\,t>B_\mu(q),\;\;\delta>0.
$$
Hence,
$$
Dim\,\overline{X}^\alpha\leq\langle\alpha,q\rangle+B_\mu(q).
$$
{\bf b. ii.} observing Lemma \ref{spectrelemme1}, assertion {\bf 2. ii.}, we obtain
$$
\mathcal{P}^{\langle\alpha,q\rangle+t+k\delta}(\underline{X}_\alpha)=0,\qquad\forall\,t>B_\mu(q),\;\;\delta>0.
$$
Hence,
$$
Dim\,\underline{X}_\alpha\leq\langle\alpha,q\rangle+t+k\delta,\qquad\forall\,t>B_\mu(q),\;\;\delta>0
$$
and finally,
$$
Dim\,\underline{X}_\alpha\leq\langle\alpha,q\rangle+B_\mu(q).
$$
\begin{lemma}\label{spectrelemme2}
$\forall\,q\in\mathbb{R}^k$ such that $\langle\alpha,q\rangle+b_\mu(q)<0$ or $\langle\alpha,q\rangle+B_\mu(q)<0$, we have $X(\alpha)=\emptyset$.
\end{lemma}
{\bf Proof.} It is based on
\begin{description}
\item Claim 1. For $q\in\mathbb{R}^k_-$ with $\langle\alpha,q\rangle+b_\mu(q)<0$ or $\langle\alpha,q\rangle+B_\mu(q)<0$, $\underline{X}_\alpha=\emptyset$.
\item Claim 2. For $q\in\mathbb{R}^k_+$ with $\langle\alpha,q\rangle+b_\mu(q)<0$ or $\langle\alpha,q\rangle+B_\mu(q)<0$, $\overline{X}^\alpha=\emptyset$.
\end{description}
Indeed, let $q\in\mathbb{R}^k_-$ and assume that $\underline{X}_\alpha\not=\emptyset$. This means that there exists at least one point $x\in\,Support(\mu)$ for which ${\underline\alpha}_{\mu_j}(x)\geq\alpha_j$, for $1\leq\,j\leq\,k$. Consequently, for all $\varepsilon>0$, there is a sequence $(r_n)_n\downarrow0$ and satisfying
$$
0<r_n<\displaystyle\frac{1}{n}\quad\hbox{and}\quad\mu_j(B(x,r_n))<r_n^{\alpha_j-\varepsilon},\;\;1\leq\,j\leq\,k.
$$
Hence,
$$
\Bigl(\mu(B(x,r_n))\Bigr)^q(2r_n)^t>2^tr_n^{\langle(\alpha-\varepsilon\1),q\rangle+t}.
$$
Choosing $t=\langle(\varepsilon\1-\alpha),q\rangle$, this induces that $\mathcal{H}_\mu^{q,t}(\{x\})>2^t$ and consequently,
$$
b_\mu(q)\geq\,dim_\mu^q(\{x\})\geq\,t,\quad\forall\,\varepsilon>0.
$$
Letting $\varepsilon\downarrow0$, it results that $b_\mu(q)\geq-\langle\alpha,q\rangle$ which is impossible. So as the first part of Claim 1. The remaining part as well as Claim 2 can be checked by similar techniques.
\begin{theorem}\label{spectretheorem1}
Let $\mu=(\mu_1,\,\mu_2,\,\dots,\,\mu_k)$ be a vector-valued Borel probability measure on $\mathbb{R}^d$ and $q\in\mathbb{R}^k$ fixed. Let further $t_q\in\mathbb{R}$, $r_q>0,\,{\underline{K}}_q,\,{\overline{K}}_q>0$, $\nu_q$ a Borel probability measure supported by $Support(\mu)$, $\varphi_q:\,\mathbb{R}_+\rightarrow\mathbb{R}$ be such that $\varphi_q(r)=o(\log\,r)$, as $r\rightarrow0$. Let finally $(r_{q,n})_n\subset]0,1[\downarrow0$ and satisfying
$$
\displaystyle\frac{\log{r_{q,n+1}}}{\log{r_{q,n}}}\rightarrow1\quad\hbox{and}\quad\displaystyle\sum_nr_{q,n}^\varepsilon<\infty,\
\forall\varepsilon>0.
$$
Assume next the following assumptions.
\begin{description}
\item[A1.] $\forall\,x\in\,Support(\mu)\,\hbox{and}\,r\in]0,r_q[$,
$$
{\underline K}_q\leq\displaystyle\frac{\nu_q(B(x,r))}{\Bigl(\mu(B(x,r))\Bigr)^q(2r)^{t_q}\exp(\varphi_q(r))}\leq{\overline K}_q.
$$
\item[A.2] $C_q(p)=\displaystyle\lim_{n\rightarrow+\infty}C_{q,n}(p)$ exists and finite for all $p\in\mathbb{R}$, where
$$
C_{q,n}(p)=\displaystyle\frac{1}{-\log{r_{q,n}}}\log\biggl(\displaystyle\int_{supp(\mu)}\Bigl(\mu(B(x,r_{q,n}))\Bigr)^pd\nu_q(x)\biggr).
$$
\end{description}
Then, the following assertions hold.
\begin{description}
\item[i.]
$$
dim(\underline{X}_{-\nabla_+C_{q}(0)}\cap\overline{X}^{-\nabla_-C_{q}(0)})\geq
$$
$$
\left\{\matrix{-\nabla_-C_{q}(0)q+\Lambda_\mu(q)\geq-\nabla_-C_q(0)q+B_\mu(q)\geq-\nabla_-C_q(0)q+b_\mu(q),\,q\in\mathbb{R}^k_-,\hfill\cr\medskip
-\nabla_+C_q(0)q+\Lambda_\mu(q)\geq-\nabla_+C_q(0)q+B_\mu(q)\geq-\nabla_+C_q(0)q+b_\mu(q),\,q\in\mathbb{R}^k_+.\hfill}\right.
$$
\item[ii.] Whenever $C_q$ is differentiable at 0, we have
$$
f_\mu(-\nabla\,C_q(0))=b_\mu^*(-\nabla\,C_q(0))=B_\mu^*(-\nabla\,C_q(0))=\Lambda_\mu^*(-\nabla\,C_q(0)).
$$
\end{description}
\end{theorem}
\begin{theorem}\label{spectretheorem2}
Assume that the hypotheses of Theorem \ref{spectretheorem1} are satisfied for all $q\in\mathbb{R}^k$. Then, the following assertions hold.
\begin{description}
\item[i.] $\alpha_\mu=-B_\mu,\quad\nu_q\,\,a.s$, whenever $B_\mu$ is differentiable at $q$.
\item[ii.] $Dom(B)\subseteq\alpha_\mu(supp(\mu))$ and $f_\mu=B_\mu^*$ on $Dom(B)$.
\end{description}
\end{theorem}
The proof of this result is based on the application of a large deviation formalism. This will permit to obtain a measure $\nu$ supported by ${\underline X}_{-\nabla_+C(0)}\cap{\overline X}^{-\nabla_-C(0)}$. To do this, we re-formulate a mixed large deviation formalism to be adapted to the mixed multifractal formalism raised in our work.  \begin{theorem}\label{largedeviation}{\bf The mixed large deviation formalism.}
Consider a sequence of vector-valued random variables $(W_n=(W_{n,1},\,W_{n,2},\,\dots,\,W_{n,k}))_n$ on a probability space $(\Omega,\,\mathcal{A},\,\mathbb{P})$ and $(a_n)_n\subset]0,+\infty[$ with $\displaystyle\lim_{n\rightarrow+\infty}a_n=+\infty$. Let next the function
$$
\matrix{C_n&:&\mathbb{R}^k\rightarrow{\overline\mathbb{R}}\hfill\cr\medskip
& &t\mapsto C_n(t)=\displaystyle\frac{1}{a_n}\log\Bigl(E(\exp(\langle\,t,W_n\rangle))\Bigr).\hfill}
$$
Assume that
\begin{description}
\item[A1.] $C_n(t)$ is finite for all $n$ and $p$.
\item[A2.] $C(t)=\displaystyle\lim_{n\rightarrow+\infty}C_n(t)$ exists and is finite for all $t$.
\end{description}
There holds that
\begin{description}
\item[i.] The function $C$ is convex.
\item[ii.] If $\nabla_-C(t)\leq\nabla_+C(t)<\alpha$, for some $t\in\mathbb{R}^k$, then
$$
\displaystyle\limsup_{n\rightarrow+\infty}\displaystyle\frac{1}{a_n}
\log\Biggl(e^{-a_nC(t)}E\biggl(\exp(\langle\,t,W_n\rangle)1_{\{\frac{W_n}{a_n}\geq\alpha\}}\biggr)\Biggr)<0.
$$
\item[iii.] If $\displaystyle\sum_ne^{-\varepsilon a_n}<\infty$ for all $\varepsilon>0$, then
$$
\displaystyle\limsup_{n\rightarrow+\infty}\frac{W_n}{a_n}\leq\nabla_+C(0)\qquad \mathbb{P}\ a.s.
$$
\item[iv.] If $\alpha<\nabla_-C(t)\leq\nabla_+C(t)$, for some $t\in\mathbb{R}^k$, then
$$
\displaystyle\limsup_{n\rightarrow+\infty}\displaystyle\frac{1}{a_n}
\log\Biggl(e^{-a_nC(t)}E\biggl(\exp(\langle\,t,W_n\rangle)1_{\{\frac{W_n}{a_n}\leq\alpha\}}\biggr)\Biggr)<0.
$$
\item[v.] If $\displaystyle\sum_ne^{-\varepsilon a_n}$ is finite for all $\varepsilon>0$, then
$$
\nabla_-C(0)\leq\displaystyle\limsup_{n\rightarrow+\infty}\frac{W_n}{a_n}\qquad\mathbb{P}\ a.s.
$$
\end{description}
\end{theorem}
\ {\bf Proof.}\\
{\bf i.} It follows from Holder's inequality.\\
{\bf ii.} Let $h\in\mathbb{R}^{*,k}_+$ be such that $C(t)+\langle\alpha,h\rangle-C(t+h)>0$. We have
$$
\matrix{&&\displaystyle\frac{1}{a_n}\log\biggl[e^{-a_nC(t)}
\mathbb{E}\Bigl(\exp(\langle\,t,W_n\rangle)1_{\{\frac{W_n}{a_n}\geq\alpha\}}\Bigr)\biggr]\hfill\cr\medskip
&=&\displaystyle\frac{1}{a_n}\log\biggl[e^{-a_nC(t)}\displaystyle\int_{\{\frac{W_n}{a_n}\geq\alpha\}}
e^{\langle\,t,W_n\rangle}d\mathbb{P}\biggr]\hfill\cr\medskip
&=&\displaystyle\frac{1}{a_n}\log\biggl[e^{-a_n(C(t)+\langle\alpha,h\rangle)}\displaystyle\int_{\{\frac{W_n}{a_n}\geq\alpha\}}
e^{\langle\,t,W_n\rangle+a_n\langle\alpha,h\rangle}d\mathbb{P}\biggr]\hfill\cr\medskip
&\leq&\displaystyle\frac{1}{a_n}\log\biggl[e^{-a_n(C(t)+\langle\alpha,h\rangle)}\displaystyle\int_{\{\frac{W_n}{a_n}\geq\alpha\}}
e^{\langle\,t+h,W_n\rangle}d\mathbb{P}\biggr]\hfill\cr\medskip
&\leq&\displaystyle\frac{1}{a_n}\log\biggl[e^{-a_n(C(t)+\langle\alpha,h\rangle)}\mathbb{E}
(\exp(\langle\,t+h,W_n\rangle))\biggr]\hfill\cr\medskip
&=&\displaystyle\frac{1}{a_n}\log\biggl[e^{-a_n(C(t)+\langle\alpha,h\rangle-C_n(t+h))}\biggr]\hfill\cr\medskip
&=&-(C(t)+\langle\alpha,h\rangle-C_n(t+h)).\hfill}
$$
Next, by taking the limsup as $n\longrightarrow+\infty$, the result follows immediately.\\

\end{document}